\documentclass[12pt,a4paper]{amsart}
\usepackage[american]{babel}
\usepackage[utf8]{inputenc}
\usepackage{amssymb,amsfonts,amsmath}
\usepackage{nicefrac,mathrsfs,mathtools}
\usepackage{enumerate}
\usepackage[pdftitle={Continuous spectral decompositions of Abelian
  group actions on C*-algebras},
  pdfauthor={Alcides Buss and Ralf Meyer},
  pdfsubject={Mathematics}
]{hyperref}
\usepackage{microtype}

\theoremstyle{plain}
\newtheorem{thm}[equation]{Theorem}
\newtheorem{prop}[equation]{Proposition}
\newtheorem{lem}[equation]{Lemma}
\newtheorem{cor}[equation]{Corollary}
\theoremstyle{definition}
\newtheorem{defn}[equation]{Definition}
\theoremstyle{remark}
\newtheorem{exmp}[equation]{Example}
\newtheorem{rem}[equation]{Remark}

\newcommand*{\Cset}{{\mathbb C}}
\newcommand*{\Nset}{{\mathbb N}}
\newcommand*{\Zset}{{\mathbb Z}}
\newcommand*{\Rset}{{\mathbb R}}

\newcommand*{\dd}{\,\mathrm{d}}

\DeclareMathOperator{\spn}{span}
\DeclareMathOperator{\cspn}{\overline{span}}
\DeclareMathOperator{\supp}{supp}

\DeclareMathOperator{\Fix}{Fix}
\DeclareMathOperator{\Prim}{Prim}
\DeclareMathOperator{\Aut}{Aut}

\newcommand*{\abs}[1]{\lvert#1\rvert}% absolute value
\newcommand*{\norm}[1]{\lVert#1\rVert}% norm
\newcommand*{\cnj}[1]{\overline{#1}}% complex conjugate
\newcommand*{\clos}[1]{\overline{#1}}% closure

\newcommand*{\ket}[1]{\lvert#1\rangle}
\newcommand*{\bra}[1]{\langle#1\rvert}
\newcommand*{\braket}[2]{\langle#1\mid#2\rangle}
\newcommand*{\KET}[1]{{\lvert#1\mathclose{\rangle\!\rangle}}}
\newcommand*{\BRA}[1]{{\mathopen{\langle\!\langle}#1\rvert}}
\newcommand*{\BRAKET}[2]{\mathopen{\langle\!\langle}#1\mid#2\mathclose{\rangle\!\rangle}}

\newcommand*{\defeq}{\mathrel{\vcentcolon=}}% defined to be equal
\newcommand*{\rc}{\overset{\mathrm{rc}}{\sim}}% relatively continuous

\newcommand*{\Cstar}{\texorpdfstring{$C^*$\nobreakdash-\hspace{0pt}}{C*-}}
\newcommand*{\Star}{\texorpdfstring{$^*$\nobreakdash-\hspace{0pt}}{*-}}
\newcommand*{\brd}{-\hspace{0pt}}
\newcommand*{\nbd}{\nobreakdash-\hspace{0pt}}

\newcommand*{\Ad}{{\mathrm{Ad}}}% adjoint action
\newcommand*{\su}{{\mathrm{su}}}% strict unconditional
\newcommand*{\uu}{{\mathrm{u}}}% unconditional
\newcommand*{\si}{{\mathrm{si}}}% square-integrable
\newcommand*{\ii}{{\mathrm{i}}}% integrable
\newcommand*{\Bound}{\mathbb B}% adjointable operators on a Hilbert module
\newcommand*{\Comp}{\mathbb K}% compact operators on a Hilbert module
\newcommand*{\Torus}{\mathbb T}% torus

\newcommand*{\CONT}{\mathcal{C}}% continuous functions
\newcommand*{\Relc}{\mathcal{R}}% relatively continuous subspace
\newcommand*{\Mult}{\mathcal{M}}% multiplier algebra
\newcommand*{\Ls}{\mathcal{L}}% used for L^p-sections of a Fell bundle
\newcommand*{\E}{\mathcal{E}}% Hilbert modules
\newcommand*{\F}{\mathcal{F}}% Hilbert modules
\newcommand*{\Hils}{\mathcal{H}}% Hilbert space
\newcommand*{\W}{\mathcal{W}}% subspace of integrable elements
\newcommand*{\IG}{\mathcal{J}(G)}%special function space on G
\newcommand*{\Bundle}[1]{\mathcal #1}%Fell bundle ...

\newcommand*{\Cross}[2]{C^*_{\mathrm{r}}(#1,#2)}% reduced crossed product

\begin{document}

\title[Continuous spectral decompositions]%
{Continuous spectral decompositions\\
  of Abelian group actions on $C^*$-algebras}

\author{Alcides Buss}
\email{abuss@math.uni-muenster.de}
\address{Mathematisches Institut\\
  Westf\"alische Wilhelms-Universit\"at M\"un\-ster\\
  Einsteinstr.\ 62\\
  48149 M\"unster\\
  Germany}

\author{Ralf Meyer}
\email{rameyer@uni-math.gwdg.de}
\address{Mathematisches Institut\\
  Georg-August-Universit\"at G\"ottingen\\
  Bunsenstr.\ 3--5\\
  37073 G\"ottingen\\
  Germany}

\begin{abstract}
  Let~\(G\) be a locally compact Abelian group and~\(\widehat{G}\) its
  Pontrjagin dual.  Following Ruy Exel, we view Fell bundles
  over~\(\widehat{G}\) as continuous spectral decompositions of
  \(G\)\nbd{}actions on \Cstar{}algebras.  We classify such spectral
  decompositions using certain dense subspaces related to Marc Rieffel's
  theory of square-integrability.

  There is a unique continuous spectral decomposition if the group acts
  properly on the primitive ideal space of the \Cstar{}algebra.  But
  there are also examples of group actions without or with several
  inequivalent spectral decompositions.
\end{abstract}


\thanks{This research was supported by the EU-Network \emph{Quantum
    Spaces and Noncommutative Geometry} (Contract HPRN-CT-2002-00280),
  the \emph{Deutsche Forschungsgemeinschaft} (SFB 478), and by CAPES.
  The first author wishes to express his thanks to Siegfried Echterhoff
  and Ruy Exel for many helpful conversations.}

\maketitle

\section{Introduction}
\label{sec:introduction}

Let~\(G\) be a locally compact Abelian group, let~\(\widehat{G}\) be its
Pontrjagin dual, and let~\(A\) be a \Cstar{}algebra with a strongly
continuous action~\(\alpha\) of~\(G\); we briefly call~\(A\) or, more
precisely, the pair \((A,\alpha)\) a \emph{\(G\)\nbd{}\Cstar{}algebra}.

Suppose first that~\(G\) is compact, and let
\begin{equation}
  \label{eq:chi_homogeneous}
  \Bundle{A}_\chi\defeq \{a\in A \mid
  \text{\(\alpha_t(a)=\braket{\chi}{t}\cdot a\)
    for all \(t\in G\)} \}
\end{equation}
for \(\chi\in\widehat{G}\), where \(\braket{\chi}{t}\defeq \chi(t)\).
For \(\chi=1\), this yields the fixed-point subalgebra
\[
\Bundle{A}_1 \defeq \{a\in A \mid
\text{\(\alpha_t(a) = a\) for all \(t\in G\)}\}.
\]
Since~\(G\) acts by \Star{}automorphisms, we have
\begin{equation}
  \label{eq:Fell_bundle_algebra}
  \Bundle{A}_\chi\cdot \Bundle{A}_\eta \subseteq \Bundle{A}_{\chi\eta}
  \quad \text{and} \quad \Bundle{A}_\chi^*=\Bundle{A}_{\chi^{-1}}
  \qquad\text{for all \(\chi,\eta\in\widehat{G}\).}
\end{equation}
Hence the spaces \(\Bundle{A}_\chi\) for \(\chi\in\widehat{G}\) form a
\emph{Fell bundle}~\(\Bundle{A}\) over~\(\widehat{G}\)
(see~\cite{fell_doran}); there is no continuity condition
because~\(\widehat{G}\) is discrete.

The cross-sectional \Cstar{}algebra \(C^*(\Bundle{A})\) of a Fell
bundle~\(\Bundle{A}\) comes with a canonical action of~\(G\) called the
\emph{dual action}~\(\widehat{\alpha}\).  We define it by
\(\widehat{\alpha}_t(a_\chi)\defeq \braket{\chi}{t}\cdot a_\chi\) for
\(t\in G\), \(\chi\in\widehat{G}\), \(a_\chi\in \Bundle{A}_\chi\); often
the convention \(\widehat{\alpha}_t(a_\chi)= \cnj{\braket{\chi}{t}}\cdot
a_\chi\) is used instead, but this does not fit
with \eqref{eq:chi_homogeneous}.

The representation theory of compact groups shows that
\(\bigoplus_{\chi\in\widehat{G}} \Bundle{A}_\chi\) is a dense subspace
of~\(A\).  Even more, the dense embedding
\(\bigoplus_{\chi\in\widehat{G}} \Bundle{A}_\chi\to A\) induces a
\(G\)\nbd{}equivariant \Star{}isomorphism \(C^*(\Bundle{A})\cong A\).

Conversely, let~\(\Bundle{A}\) be any Fell bundle over~\(\widehat{G}\)
and let~\(A\) be its cross-sectional \Cstar{}algebra, equipped with the
dual action of~\(G\).  Then the spectral decomposition
\eqref{eq:chi_homogeneous} of~\(A\) recovers the original
bundle~\(\Bundle{A}\) (up to isomorphism).  Thus \(\Bundle{A}\mapsto
C^*(\Bundle{A})\) yields an equivalence of categories between the
categories of Fell bundles over~\(\widehat{G}\) and of
\(G\)\nbd{}\Cstar{}algebras.

Now let~\(G\) be a non-compact Abelian group.  A \emph{(continuous)
  spectral decomposition} of a \(G\)\nbd{}\Cstar{}algebra~\(A\) is a
Fell bundle~\(\Bundle{B}\) over~\(\widehat{G}\) with a
\(G\)\nbd{}equivariant \Star{}isomorphism \(C^*(\Bundle{B})\cong A\).
Since~\(\widehat{G}\) is no longer discrete, a Fell bundle over~\(\widehat{G}\)
is a \emph{continuous} Banach bundle
\((\Bundle{B}_\chi)_{\chi\in\widehat{G}}\) with algebraic operations as
in \eqref{eq:Fell_bundle_algebra}; the topology on the bundle is
specified most conveniently by a space of continuous sections
(see~\cite{fell_doran}).

For non-compact~\(G\), spectral decompositions are harder to construct
because \eqref{eq:chi_homogeneous} does not help: the subspaces in
\eqref{eq:chi_homogeneous} are~\(\{0\}\) if~\(A\) is the cross-sectional
\Cstar{}algebra of a Fell bundle.  Moreover, a Fell bundle decomposition
no longer exists for general \(G\)\nbd{}\Cstar{}algebras because a
\Cstar{}algebra of the form \(C^*(\Bundle{B})\) for a Fell
bundle~\(\Bundle{B}\) over~\(\widehat{G}\) is never unital unless~\(G\) is
compact.

This gives rise to the following questions:
\begin{enumerate}[(I)]
\item Which \(G\)\nbd{}\Cstar{}algebras have a continuous spectral
  decomposition?

\item How many such decompositions are there?

\end{enumerate}

The following theorem of Ruy Exel (Theorem 11.14 in~\cite{exel1})
answers the first question.  Let \(\Mult(A)\) be the multiplier algebra
of~\(A\).

\begin{thm}
  \label{thm:main_thm_Exel}
  A \(G\)\nbd{}\Cstar{}algebra \((A,\alpha)\) has a continuous spectral
  decomposition if and only if there is a dense subspace \(\W\subseteq
  A\) with \(\W=\W^*\) and the following properties:
  \begin{description}
  \item[integrability] for all \(a\in\W\) and \(\chi\in\widehat{G}\),
    the integrals
    \[
    E_\chi(a) \defeq
    \int_G^\su\alpha_t(a) \cdot \cnj{\braket{\chi}{t}} \dd{t}
    \]
    exist unconditionally in the strict topology on \(\Mult(A)\);

  \item[relative continuity] for all \(a,b\in\W\), we have
    \[
    \lim_{\zeta\to1}
    {}\norm{E_{\chi\zeta}(a) \cdot E_\eta(b) - E_\chi(a) \cdot
      E_{\zeta\eta}(b)} = 0
    \]
    uniformly in \(\chi,\eta\in\widehat{G}\).
  \end{description}
\end{thm}

If~\(A\) is of the form \(C^*(\Bundle{B})\) for some Fell
bundle~\(\Bundle{B}\) over~\(\widehat{G}\) and~\(\W\) is the linear span
of products \(a*a^*\) for compactly supported continuous sections~\(a\),
then these conditions are satisfied by Proposition 10.2 in~\cite{exel1}.

It was not clear to Exel whether the second, technical condition in
Theorem \ref{thm:main_thm_Exel} is really necessary (compare Question
11.16 in~\cite{exel1}); its necessity was established by the second
author in~\cite{ralf1}.

The goal of this article is to classify the continuous spectral
decompositions of a \(G\)\nbd{}\Cstar{}algebra~\(A\); we show that they
correspond bijectively to dense subspaces~\(\W\) with the properties in
Theorem \ref{thm:main_thm_Exel} that are also ``complete'' in a suitable
sense.

We use this to describe all continuous spectral decompositions of
\(\Comp(L^2\widehat{G})\) with \(G\)\nbd{}action induced by the regular
representation of~\(G\); these are closely related to line bundles on
open subsets of~\(\widehat{G}\) of full measure (see
\textsection\ref{sec:ex_Hilbert}).  We will exhibit non-isomorphic Fell
bundles that provide continuous spectral decompositions of
\(\Comp(L^2\widehat{G})\).

Instead of \emph{integrable} elements as in Theorem
\ref{thm:main_thm_Exel}, we prefer to use \emph{square-integrable}
elements as in~\cite{ralf1}.  A subspace \(\Relc\subseteq A\) is
\emph{square-integrable} and \emph{relatively continuous} if and only if
\[
\W\defeq \spn \{aa^*\mid a\in\Relc\}
\]
satisfies the two conditions in Theorem \ref{thm:main_thm_Exel}.  We
call~\(\Relc\) complete if~\(\Relc\) is a \(G\)\nbd{}invariant right
ideal in~\(A\) and closed with respect to the norm
\[
\norm{\xi}_\si\defeq \norm{\xi}+\norm{E_1(\xi\xi^*)}^{1/2}.
\]
A \emph{continuously square-integrable \(G\)\nbd{}\Cstar{}algebra} is a
pair \((A,\Relc)\) consisting of a \(G\)\nbd{}\Cstar{}algebra~\(A\) and
a dense, relatively continuous, complete subspace of square-integrable
elements~\(\Relc\) \cite{ralf1}.

If~\(\Bundle{B}\) is a Fell bundle over~\(\widehat{G}\), then
\[
(A,\Relc)\defeq
\bigl(C^*(\Bundle{B}),\clos{\CONT_c(\Bundle{B})}^\si\bigr)
\]
is a continuously square-integrable \(G\)\nbd{}\Cstar{}algebra; here
\(\clos{\CONT_c(\Bundle{B})}^\si\) denotes the closure of
\(\CONT_c(\Bundle{B})\) with respect to the \(\si\)-norm.  Conversely,
given a continuously square-integrable \(G\)\nbd{}\Cstar{}algebra
\((A,\Relc)\), we use results of Exel in~\cite{exel1} to construct a
spectral decomposition \(\Bundle{B}(A,\Relc)\) of~\(A\).  Combining
techniques from \cite{exel1} and \cite{ralf1}, we show that the
isomorphism \(A\cong C^*\bigl(\Bundle{B}(A,\Relc)\bigr)\) maps~\(\Relc\)
isomorphically onto \(\clos{\CONT_c(\Bundle{B}(A,\Relc))}^\si\).

If \((A,\Relc)\) is already of the form
\(\bigl(C^*(\Bundle{B}),\clos{\CONT_c(\Bundle{B})}^\si\bigr)\) for some
Fell bundle~\(\Bundle{B}\), then we construct a canonical isomorphism
\(\Bundle{B}\cong \Bundle{B}(A,\Relc)\).  As a result, the constructions
\[
\Bundle{B}\mapsto
\bigl(C^*(\Bundle{B}),\clos{\CONT_c(\Bundle{B})}^\si\bigr),
\qquad
(A,\Relc)\mapsto \Bundle{B}(A,\Relc)
\]
are inverse to each other.  They provide an equivalence between the
categories of Fell bundles over~\(\widehat{G}\) and of continuously
square-integrable \(G\)\nbd{}\Cstar{}algebras (see Theorem
\ref{the:main}).

The main object of interest in~\cite{ralf1} is the \emph{generalized
  fixed-point algebra} of \((A,\Relc)\), which makes sense for
non-Abelian groups as well.  In our context, it agrees with the unit
fiber \(\Bundle{B_1}(A,\Relc)\) of the Fell bundle
\(\Bundle{B}(A,\Relc)\).  Thus the information needed to construct the
generalized fixed-point algebra already yields the whole spectral
decomposition.  But we cannot reconstruct~\(\Relc\) from the generalized
fixed-point algebra: there are continuous spectral decompositions of
\(\Comp(L^2\widehat{G})\) with non-isomorphic Fell bundles but with the
same generalized fixed-point algebra.

A \(G\)\nbd{}\Cstar{}algebra is called \emph{spectrally proper} if the
induced action of~\(G\) on its primitive ideal space is proper.  In this
case, it is shown in~\cite{ralf1} that there is a unique dense,
relatively continuous, complete subspace.  As a consequence, there is a
\emph{unique} continuous spectral decomposition.  We consider this case
in \textsection\ref{sec:spec_prop} and specialize even further to
\(G\)\nbd{}\Cstar{}algebras of the form \(\CONT_0(X)\) for a proper
\(G\)\nbd{}space~\(X\).  We show that this category of
\(G\)\nbd{}\Cstar{}algebras is equivalent to the category of
\emph{commutative} Fell bundles.

Marc Rieffel introduced square-integrability as a non\nbd{}commutative
generalization of proper actions on spaces \cite{rieffel2,rieffel1}.
The results in \textsection\ref{sec:spec_prop} show once more that
proper actions have special properties that are not captured by
(square)-integrability.  As already suggested in \cite{ralf2,ralf1},
square-integrability is closer to \emph{stability} than to properness.

\section{Preliminaries}
\label{sec:prelim}

In this section, we review some relevant results of
\cite{exel2,exel1,ralf2,ralf1,rieffel1} and fix our notation and
conventions.

Let~\(G\) be a locally compact group and let \((B,\beta)\) be a
\(G\)\nbd{}\Cstar{}algebra.  The completion of \(\CONT_c(G,B)\) with
respect to the \(B\)\nbd{}valued inner product
\[
\braket{f}{g}_B\defeq \int_G f(t)^*\cdot g(t) \dd{t}
\]
for \(f,g\in \CONT_c(G,B)\) is a \(G\)\nbd{}equivariant Hilbert
\(B\)\nbd{}module --- briefly \emph{Hilbert \(B,G\)-module} --- for the
\(G\)\nbd{}action
\[
(\delta_t f)(s)\defeq \beta_t\bigl(f(t^{-1}s)\bigr)
\qquad \forall\; s,t\in G, f\in \CONT_c(G,B);
\]
we denote this Hilbert \(B,G\)\brd{}module by \(L^2(G,B)\).  It is the
prototype for square-integrable Hilbert \(B,G\)-modules.

\subsection{Square-integrability}
\label{sec:si_cont_si}

Let~\((\E,\gamma)\) be a Hilbert \(B,G\)-module.  The following
constructions are taken from~\cite{ralf1}.  Given \(\xi\in\E\), we
define
\begin{alignat*}{2}
\BRA{\xi} &\colon \E\to\CONT_b(G,B),&\qquad
(\BRA{\xi}\eta)(t) &\defeq \braket{\gamma_t(\xi)}{\eta},\\
\KET{\xi} &\colon \CONT_c(G,B)\to\E,&\qquad
\KET{\xi} f &\defeq \int_G \gamma_t(\xi)\cdot f(t)\dd{t},
\end{alignat*}
We call \(\xi\in\E\) \emph{square-integrable} if \(\BRA{\xi}\eta\in
L^2(G,B)\) for all \(\eta\in\E\).  Then \(\BRA{\xi}\) becomes an
adjointable operator \(\E\to L^2(G,B)\), whose adjoint extends
\(\KET{\xi}\) to an adjointable operator \(L^2(G,B)\to\E\); we denote
these extensions by \(\BRA{\xi}\) and \(\KET{\xi}\) as well.  We also
write \(\BRAKET{\xi}{\eta}\defeq \BRA{\xi}\circ\KET{\eta}\).
Conversely, if \(\KET{\xi}\) extends to an \emph{adjointable} operator
\(L^2(G,B)\to\E\), then~\(\xi\) is square-integrable \cite{ralf1}; a
bounded extension of \(\KET{\xi}\) does not suffice for this.

\begin{rem}
  \label{rem:cont_in_Ltwo}
  What does \(\BRA{\xi}\eta\in L^2(G,B)\) mean?

  Let \((w_i)_{i\in I}\) be a net of continuous, compactly supported
  functions \(G\to[0,1]\) with \(w_i(t)\to 1\) uniformly on compact
  subsets of~\(G\).  We call \(f\in\CONT_b(G,B)\)
  \emph{square-integrable} and write \(f\in L^2(G,B)\) if the net
  \((w_i\cdot f)_{i\in I}\) in \(\CONT_c(G,B)\) is a Cauchy net in
  \(L^2(G,B)\).  We identify~\(f\) with the limit of this net in
  \(L^2(G,B)\).  It follows from Proposition
  \ref{pro:unconditional_integral} that this definition does not depend
  on the choice of the net \((w_i)_{i\in I}\) and that the continuity of
  the functions~\(w_i\) is irrelevant; it suffices if the
  functions~\(w_i\) are measurable and uniformly bounded.
\end{rem}

Let \(\E_\si\subseteq\E\) be the subspace of square-integrable elements.
We call~\(\E\) \emph{square-integrable} if~\(\E_\si\) is dense
in~\(\E\).  We equip~\(\E_\si\) and its subspaces with the norm
\[
\norm{\xi}_\si
\defeq \norm{\xi}+\norm{\KET{\xi}}
= \norm{\xi}+\norm{\BRA{\xi}}
= \norm{\braket{\xi}{\xi}}^{\nicefrac{1}{2}}
+ \norm{\BRAKET{\xi}{\xi}}^{\nicefrac{1}{2}}.
\]
This turns~\(\E_\si\) into a Banach space.  The right \(B\)\nbd{}module
structure satisfies \(\norm{\xi\cdot b}_\si\le \norm{\xi}_\si\cdot
\norm{b}\) for all \(\xi\in\E_\si\), \(b\in B\).  It is unclear
whether~\(\E_\si\) is an \emph{essential} \(B\)\nbd{}module, that is,
whether \(\E_\si\cdot B\) is dense in~\(\E_\si\).  The group~\(G\) acts
on~\(\E_\si\) by bounded operators; this action need not be continuous
and is uniformly bounded if and only if~\(G\) is unimodular.  The
convolution turns \(\CONT_c(G,B)\) into an algebra.  We furnish
\(G\)\nbd{}invariant \(B\)\nbd{}submodules of~\(\E\) such as~\(\E_\si\)
with the following \(\CONT_c(G,B)\)-module structure (see Equation (24)
in~\cite{ralf1}):
\begin{equation}
  \label{eq:CONT_GB_module}
  \xi*K \defeq \int_G \gamma_t\bigl(\xi\cdot K(t^{-1})\bigr) \dd{t}
  \qquad\forall\; \xi\in\E,\ K\in\CONT_c(G,B).
\end{equation}

\begin{defn}
  \label{def:complete_si}
  A subspace \(\Relc\subseteq\E_\si\) is called \emph{complete} if it is
  closed (for the \(\si\)-norm) and a \(\CONT_c(G,B)\)-submodule for the
  module structure defined in \eqref{eq:CONT_GB_module}.  The
  \emph{completion} of a subset \(\Relc\subseteq \E_\si\) is the
  smallest complete subset containing~\(\Relc\) or, more explicitly, the
  \(\si\)-norm closed linear span of \(\Relc\cup \Relc*\CONT_c(G,B)\).
\end{defn}

It was already noticed by Marc Rieffel that~\(\E_\si\) is too big to
construct a generalized fixed-point algebra \cite{rieffel1}.  We need
\emph{relatively continuous} subsets of~\(\E_\si\).  This notion uses
the faithful \Star{}representation
\begin{multline}
  \label{eq:cross_represent}
  \rho\colon \Cross{G}{B}\to \Bound\bigl(L^2(G,B)\bigr),
  \qquad K\mapsto \rho_K,\\
  (\rho_K f)(t)
  \defeq \int_G \beta_t\bigl(K(t^{-1}s)\bigr) \cdot f(s)\dd{s}
  \qquad\forall\; f,K\in\CONT_c(G,B),\ t\in G
\end{multline}
(see Equation (10) in~\cite{ralf1}).  Let
\(\Bound^G\bigl(L^2(G,B)\bigr)\subseteq \Bound\bigl(L^2(G,B)\bigr)\) be
the subalgebra of \(G\)\nbd{}equivariant adjointable operators on
\(L^2(G,B)\).  The range of~\(\rho\) is contained in
\(\Bound^G\bigl(L^2(G,B)\bigr)\).  We always identify \(\Cross{G}{B}\)
with its image in \(\Bound^G\bigl(L^2(G,B)\bigr)\).

\begin{defn}
  \label{def:relatively_continuous}
  A subset \(\Relc\subseteq\E_\si\) is \emph{relatively continuous} if
  \(\BRAKET{\xi}{\eta}\) belongs to \(\Cross{G}{B}\) for all
  \(\xi,\eta\in\Relc\).
\end{defn}

\begin{prop}
  \label{pro:defn_complete}
  Let~\(\E\) be a \(G\)\nbd{}equivariant Hilbert \(B,G\)\brd{}module,
  and let \(\Relc\subseteq\E_\si\) be a relatively continuous subspace.

  If\/~\(\Relc\) is complete, then the action of~\(G\) on~\(\E\)
  restricts to a continuous action on~\(\Relc\), and~\(\Relc\) is an
  essential right \(B\)\nbd{}module, that is, \(\Relc\cdot B=\Relc\).

  Conversely, \(\Relc\) is complete if it is \(\si\)-norm closed,
  \(G\)\nbd{}invariant, and a \(B\)\nbd{}submodule.
\end{prop}

\begin{proof}
  The first assertion is Proposition 6.4 in~\cite{ralf1}.  For the
  converse assertion, we notice first that the actions of \(G\)
  and~\(B\) on~\(\Relc\) are continuous because they are continuous on
  the completion of~\(\Relc\) by the first assertion.  Hence the
  function \(G\ni t\mapsto\gamma_t(\xi\cdot K(t^{-1}))\in\Relc\) for
  \(\xi\in\Relc\), \(K\in\CONT_c(G,B)\) is a
  \(\norm{\cdot}_\si\)-continuous function of compact support.
  Since~\(\Relc\) is closed, the integral \(\xi*K\) in
  \eqref{eq:CONT_GB_module} belongs to~\(\Relc\) as well; thus~\(\Relc\)
  is a closed \(\CONT_c(G,B)\)-submodule of~\(\E_\si\) as desired.%\qed
\end{proof}

\begin{cor}
  \label{cor:completion}
  Let~\(\E\) be a Hilbert \(B,G\)\nbd{}module and let
  \(\Relc\subseteq\E_\si\) be relatively continuous.  Suppose
  that~\(\Relc\) is \(G\)\nbd{}invariant and \(\Relc\cdot
  D\subseteq\clos{\Relc}^\si\) for a dense subset~\(D\) of~\(B\).  Then
  the completion of~\(\Relc\) is \(\clos{\Relc}^\si\).
\end{cor}

\begin{proof}
  The \(\si\)-norm closure \(\clos{\Relc}^\si\) is closed by
  construction, \(G\)\nbd{}invariant because~\(\Relc\) is, and a
  \(B\)\nbd{}submodule as well.  Proposition \ref{pro:defn_complete}
  shows that it is the smallest complete subspace
  containing~\(\Relc\).%\qed
\end{proof}

Given a relatively continuous subset \(\Relc\subseteq \E\), we define
\begin{equation}
  \label{eq:F_E_R}
  \F(\E,\Relc)
  \defeq \cspn{}\bigl(\KET{\Relc}\circ \Cross{G}{B}\bigr)
  \subseteq \Bound^G(L^2(G,B),\E),
\end{equation}
where \(\Bound^G(L^2(G,B),\E)\) denotes the space of equivariant,
adjointable operators \(L^2(G,B)\to\E\).  Proposition 6.1
in~\cite{ralf1} asserts that
\[
\F(\E,\Relc)\circ \Cross{G}{B} \subseteq \F(\E,\Relc),\qquad
\F(\E,\Relc)^*\circ \F(\E,\Relc) \subseteq \Cross{G}{B},
\]
so that \(\F(\E,\Relc)\) becomes a Hilbert \(\Cross{G}{B}\)-module in an
obvious way \cite{ralf1}.  Furthermore, if~\(\Relc\) is dense, then
\(\F(\E,\Relc)\) is essential in the sense that the linear span of
\(\F(\E,\Relc)\bigl(L^2(G,B)\bigr)\) is dense in~\(\E\).

We use an approximate identity in \(\Cross{G}{B}\) to show
\(\KET{\xi}\in\F(\E,\Relc)\) for all \(\xi\in\Relc\).  If~\(\Relc\) is
complete, then \(\F(\E,\Relc) = \clos{\KET{\Relc}}\).  If~\(\Relc\) is
not complete, we may replace it by its completion.

\begin{defn}
  \label{def:cont_si}
  A \emph{continuously square-integrable Hilbert \(B,G\)-module} is a
  pair \((\E,\Relc)\) consisting of a Hilbert \(B,G\)-module~\(\E\) and
  a dense subspace \(\Relc\subseteq\E\) that is contained in~\(\E_\si\)
  and complete and relatively continuous.

  The \emph{generalized fixed-point algebra} \(\Fix(\E,\Relc)\) is the
  closed linear span of \(\KET{\Relc}\BRA{\Relc}\) in \(\Bound^G(\E)\),
  the \Cstar{}algebra of \(G\)\nbd{}equivariant adjointable operators
  on~\(\E\).
\end{defn}

There is a canonical isomorphism \(\Fix(\E,\Relc) \cong
\Comp\bigl(\F(\E,\Relc)\bigr)\).

\begin{defn}
  \label{def:Relc-continuous}
  Let \((\E,\Relc)\) and \((\E',\Relc')\) be continuously
  square\brd{}integrable Hilbert \(B,G\)-modules.  An adjointable
  operator \(T\in\Bound(\E,\E')\) is called \emph{\(\Relc\)-continuous}
  if \(T(\Relc)\subseteq\Relc'\) and \(T^*(\Relc')\subseteq\Relc\).
\end{defn}

By definition, the \(\Relc\)\nbd{}continuous \(G\)\nbd{}equivariant
operators are the morphisms in the category of continuously
square-integrable Hilbert \(B,G\)-modules.  It follows from Theorem 6.2
in~\cite{ralf1} that \((\E,\Relc)\mapsto \F(\E,\Relc)\) is an
\emph{equivalence} between the categories of continuously
square-integrable Hilbert \(B,G\)\brd{}modules and of Hilbert
\(\Cross{G}{B}\)-modules.  Hence it induces a bijection between
\emph{isomorphism classes} of continuously square-integrable Hilbert
modules and Hilbert \(\Cross{G}{B}\)-modules, respectively.
% A step in the proof is that we can recover a complete relatively
% continuous subset \(\Relc\subseteq\E_\si\) from \(\F(\E,\Relc)\) by
% \[
% \Relc = \{\xi\in\E_\si\mid \KET{\xi}\in\F(\E,\Relc)\}.
% \]

We are particularly interested in the following special situation:

\begin{defn}
  \label{def:cont_si_Cstar}
  A \emph{continuously square-integrable \(G\)\nbd{}\Cstar{}algebras}
  (for a locally compact group~\(G\)) is a triple \((A,\alpha,\Relc)\),
  where \((A,\alpha)\) is a \(G\)\nbd{}\Cstar{}algebra and
  \(\Relc\subseteq A_\si\) is complete and relatively continuous and
  dense in~\(A\).  Here we view~\(A\) as a Hilbert \(A,G\)-module as
  usual.
\end{defn}

The continuously square-integrable \(G\)\nbd{}\Cstar{}algebras form a
category: morphisms \((A_1,\Relc_1)\to (A_2,\Relc_2)\) are
\(G\)\nbd{}equivariant \Star{}homomorphisms \(f\colon A_1\to A_2\) with
\(f(\Relc_1)\subseteq\Relc_2\).

\begin{rem}
  \label{rem:cont_si_Cstar_alternative}
  Since continuously square-integrable Hilbert \(B,G\)-modules
  correspond to Hilbert \(\Cross{G}{B}\)-modules, a continuously
  square\brd{}integrable \(G\)\nbd{}\Cstar{}algebra is equivalent to a
  quadruple \((A,\alpha,\F,\phi)\), where \((A,\alpha)\) is as above,
  \(\F\) is a Hilbert \(\Cross{G}{A}\)-module, and~\(\phi\) is an
  isomorphism of Hilbert \(A,G\)-modules
  \[
  \F\otimes_{\Cross{G}{A}} L^2(G,A) \cong A.
  \]
  This description is particularly useful if~\(G\) is Abelian and
  \(A=\Comp(\Hils)\) for some \(G\)\nbd{}Hilbert space~\(\Hils\) because
  then \(\Cross{G}{A}\) is Morita--Rieffel equivalent to
  \(\Cross{G}{\Cset}\cong \CONT_0(\widehat{G})\), so that
  \(\Cross{G}{A}\)-Hilbert modules are easy to classify (see
  \textsection\ref{sec:ex_Hilbert}).
\end{rem}

\subsection{Integrable versus square-integrable elements}
\label{sec:cont_si_to_Fell}

We mainly work with square-integrable elements, whereas Exel uses
integrable elements in \cite{exel2,exel1}.  Since the construction of
Fell bundles is rather technical, we do not want to redo too many of
Exel's constructions.  Instead, we use the equivalence between the
approaches via integrable and square-integrable elements, which goes
back to Rieffel \cite{rieffel1}.

\begin{defn}[\cite{exel2}]
  \label{def:loc_integrable}
  A function \(f\colon G\to A\) is \emph{locally integrable} if it is
  Bochner integrable over every measurable, relatively compact subset
  \(K\subseteq G\), that is, \(f|_K \in L^1(K,A)\).  It is
  \emph{unconditionally integrable} if it is locally integrable and the
  net \(\left(\int_K f(t)\dd{t}\right)_{K\in\mathcal{K}}\) of Bochner
  integrals converges in the norm topology of~\(A\),
  where~\(\mathcal{K}\) denotes the directed set of measurable
  relatively compact subsets of~\(G\) ordered by inclusion.  The limit
  of this net is called \emph{unconditional integral} and denoted by
  \(\int_G^\uu f(t) \dd{t}\).
\end{defn}

\begin{prop}
  \label{pro:unconditional_integral}
  Let \(f\colon G\to A\) be a locally integrable function with
  \(f(t)\geq 0\) (in~\(A\)) for all \(t\in G\).  Then the following are
  equivalent:
  \begin{itemize}
  \item \(f\) is unconditionally integrable;

  \item \(\left(\int_G w_i(t)\cdot f(t) \dd{t}\right)_{i\in I}\)
    converges in~\(A\) for all nets \((w_i)_{i\in I}\) as in Remark
    \ref{rem:cont_in_Ltwo};

  \item \(\left(\int_G w_i(t)\cdot f(t) \dd{t}\right)_{i\in I}\)
    converges in~\(A\) for some net \((w_i)_{i\in I}\) of measurable,
    compactly supported functions \(w_i\colon G\to [0,1]\) with \(
    w_i(t)\to 1\) uniformly on compact subsets of~\(G\).
  \end{itemize}
  Furthermore, in this case we have
  \[
  \lim_{i\in I} \int_G w_i(t) \cdot f(t) \dd{t}
  = \int_G^\uu f(t) \dd{t}.
  \]
\end{prop}

\begin{proof}
  The first condition implies the other conditions by Proposition 3.4
  in~\cite{exel1}.  For the converse, it is enough to prove that the
  third condition implies the first one.  Let \(a\in A_+\) be the limit
  of the net \(\left(\int_G w_i(t)\cdot f(t) \dd{t}\right)_{i\in I}\).
  Using the locally uniform convergence \(w_i\to1\), we get \(\int_K
  f(t)dt\le a\) for any measurable relatively compact subset
  \(K\subseteq G\).  Given \(i\in I\), we set \(K_i\defeq\supp(w_i)\).
  Since \(w_i\le 1\), we get
  \[
  0\le \int_G w_i(t)\cdot f(t)\dd{t}
  \le \int_{K_i} f(t)\dd{t}
  \le \int_K f(t)\dd{t}
  \le a
  \]
  for any measurable relatively compact subset~\(K\) of~\(G\) that
  contains~\(K_i\).  The assertion follows.%\qed
\end{proof}

\begin{rem}
  \label{rem:unconditional_integral}
  The first two conditions in Proposition
  \ref{pro:unconditional_integral} remain equivalent for any locally
  integrable function \(G\to A\), that is, no positivity assumption is
  needed for this.  We omit the proof.
\end{rem}

Let \((A,\alpha)\) be a \(G\)\nbd{}\Cstar{}algebra.

\begin{defn}
  \label{def:integrable}
  Let \(A_+\subseteq A\) be the subset of positive elements.  We call
  \(a\in A_+\) \emph{integrable} if the function \(G\ni
  t\mapsto\alpha_t(a)\in A\) is \emph{strictly-unconditionally
    integrable}, that is, if the functions
  \[
  G\to A,
  \qquad
  t\mapsto \alpha_t(a) \cdot b,
  \quad
  t\in G\mapsto b \cdot \alpha_t(a)
  \]
  are unconditionally integrable for each~\(b\) in~\(A\).  Let
  \(\int_G^\su\alpha_t(a)\dd{t}\) be the multiplier \((L,R)\) of~\(A\)
  given by
  \[
  L(b)=\int_G^\uu \alpha_t(a)\cdot b \dd{t},
  \qquad
  R(b)=\int_G^\uu b\cdot \alpha_t(a) \dd{t}.
  \]
  Let \(A_\ii\subseteq A\) be the subspace spanned by the integrable
  elements in~\(A_+\).
\end{defn}

\begin{rem}
  \label{lem:positive_integrable}
  A positive element belongs to~\(A_\ii\) if and only if it is
  integrable.
\end{rem}

Let \((\E,\gamma)\) be a Hilbert \(B,G\)\brd{}module.  Let \(A \defeq
\Comp(\E)\) be the \(G\)\nbd{}\Cstar{}algebra of compact operators
on~\(\E\) with~\(G\) acting by \(\alpha_t(a) \defeq \gamma_t\circ
a\circ\gamma_t^{-1}\) for all \(a\in \Comp(\E)\) and \(t\in G\).  For
\(\xi,\eta\in\E\), we define a compact operator
\[
\ket{\xi}\bra{\eta}\colon \E\to\E,
\qquad
\zeta\mapsto \xi\cdot\braket{\eta}{\zeta}_B.
\]
We have \(\alpha_t(\ket{\xi}\bra{\eta}) =
\ket{\gamma_t\xi}\bra{\gamma_t\eta}\) for all \(\xi,\eta\in\E\).  Let
\[
M\colon L^\infty(G)\to \Bound\bigl(L^2(G,B)\bigr),\qquad
\phi\mapsto M_\phi
\]
be the representation by pointwise multiplication operators.

\begin{prop}
  \label{pro:Comp_integrable}
  Let~\(\E\) be a Hilbert \(B,G\)\nbd{}module.
  \begin{enumerate}[(i)]
  \item If \(\xi,\eta\in\E_\si\) and \(\phi\in L^\infty(G)\), then
    \(\ket{\xi}\bra{\eta}\in\Comp(\E)_\ii\) and
    \[
    \int_G^\su \phi(t) \cdot \alpha_t(\ket{\xi}\bra{\eta}) \dd{t}
    = \KET{\xi} \circ M_\phi \circ \BRA{\eta}
    \qquad \text{in \(\Bound(\E)\).}
    \]

  \item Conversely, if \(\ket{\xi}\bra{\xi}\in\Comp(\E)_\ii\) then
    \(\xi\in \E_\si\).
  \end{enumerate}
  Thus \(\Relc\subseteq\E_\si\) if and only if \(\ket{\Relc}\bra{\Relc}
  \defeq \{\ket{\xi}\bra{\eta}\mid \xi,\eta\in\Relc\} \subseteq
  \Comp(\E)_\ii\).
\end{prop}

\begin{proof}
  Polarization allows us to assume that \(\xi=\eta\) and \(\phi\ge0\).
  Let \((w_i)_{i\in I}\) be a net as in Remark \ref{rem:cont_in_Ltwo}.
  Consider the net
  \[
  I_i\defeq \int_G w_i(t) \cdot \alpha_t(\ket{\xi}\bra{\eta}) \dd{t}
  \]
  in \(\Comp(\E)\).  For \(\zeta\in\E\), we compute
  \begin{multline*}
    I_i(\zeta)
    = \int_G w_i(t) \cdot \gamma_t(\xi)
    \braket{\gamma_t(\eta)}{\zeta} \dd{t}
    = \int_G \gamma_t(\xi) \cdot  w_i(t) \cdot
    \BRA{\eta}(\zeta)(t) \dd{t}
    \\ = \KET{\xi} \bigl(w_i\cdot\BRA{\eta}(\zeta)\bigr)
    = \KET{\xi}\circ M_{w_i} \circ \BRA{\eta}(\zeta).
  \end{multline*}
  Hence \((I_i)_{i\in I}\) is a bounded net with \(\norm{I_i}\le
  \norm{\KET{\xi}} \cdot \norm{\BRA{\eta}}\); it converges pointwise to
  \(\KET{\xi}\BRA{\eta}\) because
  \(\KET{\xi}\bigl(w_i\cdot\BRA{\eta}(\zeta)\bigr)\) converges to
  \((\KET{\xi}\BRA{\eta})(\zeta)\) by Remark \ref{rem:cont_in_Ltwo}.  It
  follows that \((I_i)_{i\in I}\) converges in the strict topology of
  \(\Comp(\E)\) to \(\KET{\xi}\BRA{\eta}\).  Since \(\xi=\eta\), we have
  \(\ket{\xi}\bra{\eta}\ge0\), so that Proposition
  \ref{pro:unconditional_integral} applies.  We get
  \(\ket{\xi}\bra{\eta}\in\Comp(\E)_\ii\) and
  \[
  \int_G^\su \alpha_t(\ket{\xi} \bra{\eta}) \dd{t}
  = \KET{\xi}\BRA{\eta}
  \]
  as desired.  The statement with \(\phi\in L^{\infty}(G)\) is proved
  similarly.  This finishes the proof of~(i).  Assertion~(ii) is Lemma
  8.1.(v) in~\cite{ralf2}.%\qed
\end{proof}

Let \((A,\alpha)\) be a \(G\)\nbd{}\Cstar{}algebra viewed as a Hilbert
module over itself.  Then the map \(\ket{a}\bra{b}\mapsto a\cdot b^*\)
induces an isomorphism \(\Comp(\E)\cong A\).  Hence Proposition
\ref{pro:Comp_integrable} and polarization yield
\[
\Relc\subseteq A_\si
\quad\iff\quad
\{a\cdot a^*\mid a\in\Relc\}\subseteq A_\ii.
\]

\subsection{Relative continuity via Fourier coefficients}
\label{sec:relc_Fourier}

>From now on, we suppose that~\(G\) is Abelian, and we
let~\(\widehat{G}\) be its dual group.  We write
\(\braket{\chi}{t}\defeq \chi(t)\) for \(\chi\in\widehat{G}\), \(t\in
G\) to emphasize the symmetry between \(G\) and~\(\widehat{G}\).  The
following results are due to Ruy Exel \cite{exel1}.

The literature contains several conventions for the dual action and the
Fourier transform, which differ by the inversion automorphism
\(\chi\mapsto\chi^{-1} = \cnj{\chi}\) on~\(\widehat{G}\).  We follow the
conventions of \cite{exel2}, which are different from those of
\cite{exel1}.  Thus the \emph{Fourier transform} of a function \(f\colon
G\to\Cset\) is defined by
\begin{equation}
  \label{eq:Fourier_scalar}
  \hat{f}\colon \widehat{G}\to\Cset,
  \qquad\hat{f}(\chi) \defeq
  \int_G f(t) \cdot \cnj{\braket{\chi}{t}}\dd{t}.
\end{equation}

\begin{defn}[Definition 6.1 of~\cite{exel1}]
  \label{def:Fourier_coeff}
  Let~\(A\) be a \(G\)\nbd{}\Cstar{}algebra.  For \(a\in A_\ii\) and
  \(\chi\in\widehat{G}\), we define
  \[
  \hat{a}(\chi) = E_\chi(a) \defeq
  \int_G^\su \alpha_t(a) \cdot \cnj{\braket{\chi}{t}} \dd{t}
  \in \Mult(A).
  \]
\end{defn}

The resulting function \(\hat{a}\colon \widehat{G}\to\Mult(A)\) is
bounded and uniformly continuous for the strict topology on \(\Mult(A)\)
by Propositions 6.2--3 of~\cite{exel1}.

We denote the unique extension of~\(\alpha\) to \(\Mult(A)\) again
by~\(\alpha\); this extension need not be continuous any more.  The
\emph{\(\chi\)\nbd{}spectral subspace} \(\Mult_\chi(A)\) of \(\Mult(A)\)
is defined by
\[
\Mult_\chi(A)\defeq \{m\in\Mult(A)\mid
\text{\(\alpha_t(m)=\braket{\chi}{t}\cdot m\) for all \(t\in G\)}
\}
\]
as in \eqref{eq:chi_homogeneous}.  It is easy to see that \(E_\chi(a)
\in \Mult_\chi(A)\), that is,
\begin{equation}
  \label{eq:x-spectral_sub}
  \alpha_t\bigl(E_\chi(a)\bigr)
  = \braket{\chi}{t}\cdot E_\chi(a)
  \quad \text{for all \(t\in G\), \(\chi\in \widehat{G}\),
    \(a\in A_\ii\).}
\end{equation}

\begin{defn}[Definition 8.1 of~\cite{exel1}]
  \label{def:relatively_continuous_int}
  Let~\(A\) be a \(G\)\nbd{}\Cstar{}algebra.  A pair of integrable
  elements \((a,b)\) is \emph{relatively continuous} if
  \[
  \lim_{\eta\to1} \norm{E_{\chi\eta}(a) \cdot E_\zeta(b) -
    E_\chi(a) \cdot E_{\eta\zeta}(b)}
  = 0
  \qquad \text{uniformly in \(\chi,\zeta\in\widehat{G}\).}
  \]
  We write \(a\rc b\) in this case.  A subset \(\W\subseteq A_\ii\) is
  called \emph{relatively continuous} if \(a\rc b\) for all \(a,b\in
  \W\).
\end{defn}

Definitions \ref{def:relatively_continuous} and
\ref{def:relatively_continuous_int} define relative continuity for
square-integrable and integrable elements in terms of \(G\)
and~\(\widehat{G}\), respectively.  We need some preparations to relate
these two notions.

Let~\(B\) be a \(G\)\nbd{}\Cstar{}algebra and let
\(\chi\in\widehat{G}\).  We let~\(M_\chi\) be the pointwise
multiplication operator
\[
(M_\chi f)(t) \defeq \cnj{\braket{\chi}{t}}\cdot f(t)
\qquad \text{for \(f\in L^2(G,B)\), \(t\in G\).}
\]
This defines a strongly continuous unitary representation \(\chi\mapsto
M_\chi\) of~\(\widehat{G}\) on \(L^2(G,B)\); this operator is
denoted~\(V_\chi\) in~\cite{exel1}.

Recall that the \emph{dual action}~\(\widehat{\beta}\) of~\(\widehat{G}\) on
the crossed product \(\Cross{G}{B}\) is defined by continuous extension
of the action \((\widehat{\beta}_\chi f)(t) \defeq \braket{\chi}{t}\cdot
f(t)\) for \(f\in\CONT_c(G,B)\), \(t\in G\), \(\chi\in\widehat{G}\).
The representation~\(\rho\) in \eqref{eq:cross_represent}
relates~\(\widehat{\beta}\) and the operators~\(M_\chi\):
\begin{equation}
  \label{eq:Mx_rhoK_Mx}
  M_\chi \circ \rho_K \circ M_\chi^* =
  \rho_{\widehat{\beta}_\chi(K)}
  \qquad \forall\; \chi\in\widehat{G},\ K\in\Cross{G}{B};
\end{equation}
it suffices to check this for \(K\in\CONT_c(G,B)\).  In particular, we
get
\begin{equation}
  \label{eq:Mx_T_Mx}
  M_\chi \circ T \circ M_\chi^*\in \Cross{G}{B}
  \quad \forall\; T\in\Cross{G}{B}.
\end{equation}

Let~\(\E\) be a Hilbert \(B,G\)\brd{}module for a
\(G\)\nbd{}\Cstar{}algebra~\(B\).  We apply the above definitions and
results to \(A=\Comp(\E)\), so that \(\Mult(A) = \Bound(\E)\) with~\(G\)
acting by conjugation: \(\alpha_t(T)=\gamma_t\circ T\circ\gamma_t^{-1}\)
for all \(T\in\Bound(\E)\), \(t\in G\).  Proposition
\ref{pro:Comp_integrable} and Equation \eqref{eq:x-spectral_sub} yield
\begin{align}
  \label{eq:E_x_ketbra}
  E_\chi(\ket{\xi}\bra{\eta})
  &= \KET{\xi} \circ M_\chi \circ \BRA{\eta}
  \\
  \label{eq:act_on_Echi}
  \gamma_t \circ E_\chi(\ket{\xi}\bra{\eta}) \circ \gamma_t^{-1}
  &= \braket{\chi}{t} \cdot E_\chi(\ket{\xi}\bra{\eta})
\end{align}
for all \(\chi\in\widehat{G}\), \(\xi,\eta\in\E_\si\), \(t\in G\)
(compare Lemma 7.4 in~\cite{exel1}).  The following proposition relates
our two notions of relative continuity in terms of \(G\)
and~\(\widehat{G}\) and is one of the main technical results
in~\cite{exel1}.

\begin{prop}
  \label{pro:compare_si_int}
  Let~\(G\) be Abelian and let~\(\E\) be a Hilbert \(B,G\)\brd{}module
  for a \(G\)\nbd{}\Cstar{}algebra~\(B\).  If \(\xi,\eta\in\E_\si\),
  then
  \[
  \ket{\xi}\bra{\xi} \rc \ket{\eta}\bra{\eta} \quad\iff\quad
  \BRAKET{\xi}{\eta}\in\Cross{G}{B}
  \]
  A set of square-integrable elements \(\Relc\subseteq\E_\si\) is
  relatively continuous (Definition \ref{def:relatively_continuous}) if
  and only if the set of integrable elements
  \begin{equation}
    \label{eq:W_Relc}
    \W_\Relc \defeq \{\ket{\xi}\bra{\eta}\mid \xi,\eta\in \Relc\}
    \subseteq\Comp(\E)_\ii
  \end{equation}
  is relatively continuous
  (Definition \ref{def:relatively_continuous_int}).
\end{prop}

\begin{proof}
  The first assertion follows as in the proof of Theorem 7.5
  in~\cite{exel1} with appropriate changes of notation, using
  \eqref{eq:E_x_ketbra}.  It implies the second assertion because the
  relative continuity of~\(\W_\Relc\) is equivalent to the relative
  continuity of \(\{\ket{\xi}\bra{\xi}\mid \xi\in\Relc\}\) by
  polarization.%\qed
\end{proof}

\subsection{Spectral invariance}
\label{sec:spectrally_invariant}

As before, the group~\(G\) is Abelian with dual group~\(\widehat{G}\),
and \(A=\Comp(\E)\) for a Hilbert \(B,G\)\brd{}module \((\E,\gamma)\)
and a \(G\)\nbd{}\Cstar{}algebra \((B,\beta)\), with \(\alpha_t(a) =
\gamma_t\circ a\circ \gamma_t^{-1}\) for all \(t\in G\), \(a\in A\).

\begin{lem}
  \label{lem:X_preserves_si}
  Given \(\xi,\eta \in \E_\si\) and \(\chi\in\widehat{G}\), let \(X
  \defeq E_\chi(\ket{\xi}\bra{\eta})\in\Bound(\E)\).  Then
  \(X(\E_\si)\subseteq \E_\si\) and
  \begin{equation}
    \label{eq:compute_X}
    \KET{X(\zeta)} = X \circ \KET{\zeta} \circ M_\chi^*
    = \KET{\xi} \circ M_\chi \circ \BRAKET{\eta}{\zeta} \circ M_\chi^*
  \end{equation}
  for all \(\zeta\in\E_\si\).
\end{lem}

\begin{proof}
  Equation \eqref{eq:act_on_Echi} yields \(\gamma_t\bigl(X(\zeta)\bigr)
  = \braket{\chi}{t} \cdot X\bigl(\gamma_t(\zeta)\bigr)\).  Now we
  verify \eqref{eq:compute_X} on \(f\in\CONT_c(G,B)\):
  \begin{multline*}
    \KET{X(\zeta)} (f)
    =\int_G \gamma_t\bigl(X(\zeta)\bigr)\cdot f(t) \dd{t}
    = \int_G \braket{\chi}{t} \cdot X\bigl(\gamma_t(\zeta)\bigr)
    \cdot f(t) \dd{t}\\
    = X\left(\int_G \gamma_t(\zeta) \cdot (M_\chi^* f)(t) \dd{t}\right)
    = X\circ \KET{\zeta}\circ M_\chi^*(f).
  \end{multline*}
  Since \(\zeta\in\E_\si\) and \(M_\chi\in\Bound\bigl(L^2(G,B)\bigr)\),
  this extends to an adjointable operator \(L^2(G,B)\to \E\).  Hence
  \(X(\zeta)\in\E_\si\) and the first equality in \eqref{eq:compute_X}
  holds.  The second one follows from \eqref{eq:E_x_ketbra}.%\qed
\end{proof}

\begin{defn}
  \label{def:spectrally_invariant}
  A subset \(\W\subseteq A_\ii\) is \emph{spectrally invariant} if
  \(E_\chi(a) \cdot b\) and \(b\cdot E_\chi(a)\) belong to~\(\W\) for
  all \(a,b\in\W\), \(\chi\in\widehat{G}\).
\end{defn}

\begin{prop}
  \label{pro:X_R-continuous}
  Let \(\Relc\subseteq\E_\si\) be relatively continuous and complete,
  and define~\(\W_\Relc\) as in \eqref{eq:W_Relc}.  Let
  \(\xi,\eta\in\Relc\), \(\chi\in\widehat{G}\), and define~\(X\) as in
  Lemma \ref{lem:X_preserves_si}.  Then \(X(\Relc)\subseteq\Relc\) and
  \(X^*(\Relc)\subseteq\Relc\).  Hence~\(X\) is
  \(\Relc\)\nbd{}continuous and~\(\W_\Relc\) is spectrally invariant.
\end{prop}

\begin{proof}
  Recall that \(T\circ \ket{\xi}\bra{\eta} = \ket{T(\xi)}\bra{\eta}\)
  and \(\ket{\xi}\bra{\eta} \circ T = \ket{\xi}\bra{T^*(\eta)}\) for all
  \(T\in\Bound(\E)\).  The spectral invariance of~\(\W_\Relc\) means
  that
  \begin{align*}
    E_\chi(\ket{\xi}\bra{\eta})\circ \ket{\zeta}\bra{\upsilon}
    &= \ket{E_\chi(\ket{\xi}\bra{\eta})(\zeta)}\bra{\upsilon},\\
    \ket{\zeta}\bra{\upsilon} \circ E_\chi(\ket{\xi}\bra{\eta})
    &= \ket{\zeta} \bra{E_\chi(\ket{\xi}\bra{\eta})^*(\upsilon)}
  \end{align*}
  belong to~\(\W_\Relc\) for all \(\xi,\eta,\zeta,\upsilon\in\Relc\).
  Equivalently, \(X(\Relc)\subseteq\Relc\) and
  \(X^*(\Relc)\subseteq\Relc\).  Since \(X^*= \KET{\eta}\circ
  M_\chi^*\circ\BRA{\xi} = E_{\cnj{\chi}}(\ket{\eta}\bra{\xi})\) has the
  same form, it suffices to prove \(X(\Relc)\subseteq\Relc\).  Let
  \(\F\defeq \F(\E,\Relc)\).  Proposition 6.3 in~\cite{ralf1} yields
  \(\Relc = \Relc_\F\).  Hence we must show \(\KET{X(\zeta)}\in\F\) for
  all \(\zeta\in\Relc\).

  Equation \eqref{eq:compute_X} yields
  \[
  \KET{X(\zeta)}
  = X \circ \KET{\zeta} \circ M_\chi^*
  = \KET{\xi} \circ M_\chi \circ \BRAKET{\eta}{\zeta} \circ M_\chi^*.
  \]
  Since~\(\Relc\) is relatively continuous, we have
  \(\BRAKET{\eta}{\zeta}\in \Cross{G}{B}\).  Now \eqref{eq:Mx_T_Mx}
  yields \(M_\chi \circ \BRAKET{\eta}{\zeta}\circ M_\chi^*\in
  \Cross{G}{B}\).  Since \(\F\circ\Cross{G}{B}\subseteq\F\), we get
  \(\KET{X(\zeta)}\in\F\) as desired.%\qed
\end{proof}

\subsection{A Fourier inversion theorem}
\label{sec:Fourier_inversion}

Let \(a\in A_\ii\) and recall that
\[
\hat{a}(\chi)= E_\chi(a)
\defeq \int_G^\su \alpha_t(a) \cdot \cnj{\braket{\chi}{t}} \dd{t}
\]
The usual Fourier inversion theorem leads us to expect
\begin{equation}
  \label{eq:Fourier-inversion_E}
  \int_{\widehat{G}} \hat{a}(\chi) \cdot \braket{\chi}{t} \dd{\chi}
  = \alpha_t(a)
  \qquad \forall\; t\in G
\end{equation}
(for suitably normalized Haar measures on \(G\) and~\(\widehat{G}\)).
In fact, Proposition 6.6 in~\cite{exel1} asserts this whenever the
integral on the left hand side converges absolutely, that is, if \(a\in
A_\ii\) and \(\int_{\widehat{G}} {}\norm{\hat{a}(\chi)} \dd{\chi} <
\infty\); this ensures existence of the integral because~\(\hat{a}\) is
strictly continuous.

We need some preparation to find enough elements \(a\in A_\ii\) with
absolutely integrable Fourier transform~\(\hat{a}\).  We
specialise \eqref{eq:CONT_GB_module} to turn~\(A\) into a right
\(L^1(G)\)-module by
\[
a*f \defeq \int_G \alpha_t(a)\cdot f(t^{-1}) \dd{t}
\qquad \forall\; a\in A, f\in L^1(G).
\]

\begin{prop}
  \label{pro:fourier_and_products}
  We have \(a*f\in A_\ii\) and
  \[
  \widehat{a*f}=\hat{a}\cdot\hat{f} \qquad \text{(pointwise product)}
  \]
  for all \(a\in A_\ii\) and \(f\in L^1(G)\), with~\(\hat{f}\) as
  in \eqref{eq:Fourier_scalar}.
\end{prop}

\begin{proof}
  Formally, this follows by a change of variables and Fubini's Theorem:
  \begin{multline*}
    \widehat{a*f}(\chi)
    = \int_G \alpha_t(a*f) \cdot \cnj{\braket{\chi}{t}}\dd{t}
    = \int_{G\times G} \alpha_{ts}(a)\cdot
    f(s^{-1}) \cdot \cnj{\braket{\chi}{t}}\dd{s}\dd{t}
    \\=\int_{G\times G} \alpha_t(a) \cdot \cnj{\braket{\chi}{t}} \cdot
    f(s^{-1}) \cdot \braket{\chi}{s}\dd{s}\dd{t}
    = \hat{a}(\chi) \cdot \hat{f}(\chi).
  \end{multline*}
  This computation requires a justification because some of the relevant
  integrals only exist in a rather weak sense.  For \(f\in\CONT_c(G)\),
  this is done in the proof of Proposition 6.5 in~\cite{exel1}.  The
  only point in Exel's argument that requires additional care occurs in
  the proof of Proposition 3.6 in~\cite{exel1}: we must show that
  \((1_K*f)_{K\in\mathcal{K}}\) is a bounded net in \(L^\infty(G)\) that
  converges uniformly on compact subsets to the constant function taking
  the value \(\int_G f(s^{-1}) \dd{s} = \int_G f(s) \dd{s}\);
  here~\(1_K\) denotes the characteristic function of~\(K\).  Since
  \[
  (1_K*f)(t)
  = \int_G 1_K(ts^{-1}) f(s) \dd{s}
  = \int_{K^{-1}t} f(s)\dd{s},
  \]
  our net is bounded with \(\norm{1_K*f}_\infty\le \norm{f}_1\) for all
  \(K\in\mathcal{K}\).  It also follows that
  \[
  \left\lvert (1_K*f)(t) - \int_G f(s) \dd{s} \right\rvert \le
  \int_{G\backslash K^{-1}t} \abs{f(s)} \dd{s}.
  \]
  For any \(\varepsilon>0\) there is a compact subset \(K_0\subseteq G\)
  with \(\int_{G\backslash K_0} \abs{f(s)} \dd{s}<\varepsilon\) because
  \(f\in L^1(G)\).  Let \(L\subseteq G\) be compact.  If \(K\subseteq
  G\) is a measurable relatively compact subset with \(L\cdot
  K_0^{-1}\subseteq K\), then \(K_0\subseteq K^{-1}t\) for all \(t\in
  L\), so that
  \[
  \int_{G\backslash K^{-1}t} \abs{f(s)} \dd{s}<\varepsilon
  \]
  for all \(t\in L\).  This means that the net
  \((1_K*f)_{K\in\mathcal{K}}\) converges to \(\int_G f(s) \dd{s}\)
  uniformly on~\(L\).%\qed
\end{proof}

\begin{defn}
  \label{def:IG}
  Let \(\IG\defeq \{f\in L^1(G)\mid \hat{f}\in \CONT_c(\widehat{G})\}\).
\end{defn}

\begin{lem}
  \label{lem:IG}
  The space \(\IG\) is dense in \(L^1(G)\), and
  \[
  \widehat{\IG}\defeq \{\hat{f}\mid f\in\IG\}
  = \{g\in \CONT_c(\widehat{G}) \mid \check{g}\in L^1(G)\}
  \]
  is dense in \(\CONT_c(\widehat{G})\) in the inductive limit topology;
  here~\(\check{g}\) denotes the inverse Fourier transform of~\(g\).
\end{lem}

\begin{proof}
  The space~\(\IG\) contains \(h=\check{f}_1\cdot \check{f}_2\) for all
  \(f_1,f_2\in\CONT_c(\widehat{G})\) because \(\hat{h} = f_1*f_2\) is
  still compactly supported and
  \[
  f_1,f_2\in L^2(\widehat{G}) \quad \Longrightarrow \quad
  \check{f}_1,\check{f}_2\in L^2(G) \quad \Longrightarrow\quad
  h = \check{f}_1 \cdot \check{f}_2 \in L^1(G).
  \]
  The products \(\check{f_1}\cdot \check{f_2}\) with
  \(f_1,f_2\in\CONT_c(\widehat{G})\) form a dense subset of \(L^1(G)\)
  because the inverse Fourier transform of \(\CONT_c(\widehat{G})\) is
  dense in \(L^2(G)\).  Similarly, \(\widehat{\IG}\) contains
  \(\CONT_c(\widehat{G})*\CONT_c(\widehat{G})\), which is dense in
  \(\CONT_c(\widehat{G})\).%\qed
\end{proof}

\begin{prop}
  \label{pro:Fourier-inversion}
  The Fourier inversion formula \eqref{eq:Fourier-inversion_E} applies
  to \(a*f\) for all \(a\in A_\ii\), \(f\in\IG\).
\end{prop}

\begin{proof}
  Proposition \ref{pro:fourier_and_products} shows that \(a*f\) has
  compactly supported Fourier transform for all \(f\in\IG\).  Hence the
  assertion follows from Proposition 6.6 in~\cite{exel1}.%\qed
\end{proof}

\begin{rem}
  \label{rem:Rieffel_integrable_Fourier_inversion}
  It is well-known that we can find an approximate identity~\((u_n)\) in
  \(L^1(G)\) such that~\(\hat{u}_n\) has compact support and \(0\le
  \hat{u}_n\le 1\) for all \(n\in\Nset\), and \(\lim \hat{u}_n=1\)
  locally uniformly on~\(\widehat{G}\).  Since the Fourier inversion
  formula \eqref{eq:Fourier-inversion_E} holds for \(a*u_n\) for all
  \(n\in\Nset\), we get
  \[
  \alpha_t(a) = \lim_{n\to\infty} u_n*\alpha_t(a) = \lim_{n\to\infty}
  \int_{\widehat{G}} \hat{a}(\chi) \cdot \hat{u}_n(\chi) \cdot
  \braket{\chi}{t} \dd{\chi}.
  \]
  Hence the Fourier inversion formula always holds if we interpret the
  integral suitably.  The only issue is the quality of the convergence.
\end{rem}

\section{Continuous spectral decompositions}
\label{csi_versus_spectral_decomp}

Now we discuss how to construct continuous spectral decompositions out
of dense, complete, relatively continuous subspaces, and vice versa.
Much of this goes back to Ruy Exel \cite{exel2,exel1}.  Throughout this
section, \(G\) is an Abelian locally compact group and~\(\widehat{G}\)
is its dual group.

\subsection{Construction of continuous spectral decompositions}
\label{sec:construct_Fell}

\begin{defn}
  \label{def:cont_spec_decomp}
  Let~\(A\) be a \(G\)\nbd{}\Cstar{}algebra.  A \emph{continuous
    spectral decomposition} of~\(A\) is a Fell bundle~\(\Bundle{B}\)
  with a \(G\)\nbd{}equivariant \Cstar{}algebra isomorphism
  \(C^*(\Bundle{B})\to A\).
\end{defn}

A similar notion already appears in Definition VIII.17.4 in
\cite{fell_doran}, where it is called a \Cstar{}bundle structure
for~\(A\) over~\(\widehat{G}\); but in \cite{fell_doran}, \(A\) does not
carry an action of~\(G\) to begin with, so that the interpretation as a
spectral decomposition is missing.

Recall that \(C^*(\Bundle{B})\) is the \Cstar{}envelope of the Banach
\Star{}algebra \(\Ls^1(\Bundle{B})\) of measurable sections~\(b\) with
\(\int_{\widehat{G}} {}\norm{b(\chi)} \dd{\chi}<\infty\).  We equip
\(C^*(\Bundle{B})\) with the \emph{dual action}~\(\beta\) of~\(G\)
defined by
\[
(\widehat{\beta}_t f)(\chi) = \braket{\chi}{t} \cdot f(\chi)
\qquad
\forall\; f\in\Ls^1(\Bundle{B}), t\in G, \chi\in\widehat{G}.
\]
We often use the subspace \(\CONT_c(\Bundle{B})\) of continuous sections
with compact support.  This is a dense, \(G\)\nbd{}invariant
\Star{}subalgebra in \(\Ls^1(\Bundle{B})\) and \(C^*(\Bundle{B})\).

Unless~\(G\) is compact, a continuous spectral decomposition cannot
exist for a trivial action because dual actions are never trivial.  Thus
we need more structure to construct a continuous spectral decomposition
for a \(G\)\nbd{}\Cstar{}algebra~\(A\).  Ruy Exel uses a subset~\(\W\)
with the following properties:
\begin{enumerate}[(i)]
\item \(\W\) is a dense linear subspace of~\(A\) with \(\W=\W^*\);

\item \(\W\subseteq A_\ii\);

\item \(\W\) is relatively continuous;

\item \(\W\) is spectrally invariant.
\end{enumerate}
This yields a continuous spectral decomposition
\(C^*\bigl(\Bundle{B}(\W)\bigr)\) for~\(A\) (see \textsection11
of~\cite{exel1}); the density of~\(\W\) in~\(A\) ensures that the map
\(C^*\bigl(\Bundle{B}(\W)\bigr)\to A\) is an isomorphism.

The fibers of \(\Bundle{B}(\W)\) are
\[
\Bundle{B}_\chi(\W)
\defeq \cspn {}\{\hat{a}(\chi)\mid a\in\W\}
\subseteq \Mult_\chi(A) \subseteq  \Mult(A)
\]
for \(\chi\in\widehat{G}\); more precisely, we take the norm closure
here.  The algebraic operations are inherited from~\(\Mult(A)\); the
topology is defined so that the sections~\(\hat{a}\) for \(a\in\W\)
\emph{generate} the space of continuous sections.  We must make this
more precise (see also II.13 in \cite{fell_doran}).

To begin with, a Fell bundle~\(\Bundle{B}\) is determined uniquely by
the space \(\CONT_c(\Bundle{B})\) of its continuous, compactly supported
sections (since continuity is a local issue, the support condition is
irrelevant here).  The space of continuous sections of~\(\Bundle{B}\) is
closed under pointwise addition and multiplication by continuous
functions, that is, it is a \(\CONT(\widehat{G})\)-module.  If
\(U\subseteq\widehat{G}\) is a relatively compact, open subset, we let
\(\CONT_0(U,\Bundle{B})\subseteq \CONT_c(\Bundle{B})\) be the subspace
of continuous sections that vanish outside~\(U\).  These subspaces are
complete with respect to the norm
\[
\norm{f}^U_\infty \defeq \sup_{\chi\in U} {}\norm{f(\chi)}.
\]
We have \(\CONT_c(\Bundle{B}) = \bigcup \CONT_0(U,\Bundle{B})\) and
equip \(\CONT_c(\Bundle{B})\) with the inductive limit topology.  These
norms and the inductive limit topology also make sense on the
corresponding larger spaces \(D_0(U,\Bundle{B})\) and
\(D_c(\Bundle{B})\) of bounded sections of~\(\Bundle{B}\).

\begin{lem}
  \label{lem:Fell_cont_section}
  The subspace \(\CONT_c(\Bundle{B}) \subseteq D_c(\Bundle{B})\) is the
  closed linear span of the set of all \(g\cdot \hat{a}\) with
  \(g\in\CONT_c(\widehat{G})\) and \(a\in\W\).  We may further restrict
  to \(\hat{g}\cdot\hat{a}\) with \(g\in\IG\) and \(a\in\W\), with
  \(\IG\) as in Definition \ref{def:IG}.
\end{lem}

\begin{proof}
  Recall that \(\CONT_c(\Bundle{B})\) is a closed
  \(\CONT_c(\widehat{G})\)-submodule of \(D_c(\Bundle{B})\) for any Fell
  bundle~\(\Bundle{B}\).  It remains to check that any continuous
  compactly supported section of~\(\Bundle{B}\) can be approximated
  uniformly by linear combinations of pointwise products \(g\cdot
  \hat{a}\) with \(g\in\CONT_c(\widehat{G})\) and \(a\in\W\).  By
  definition, a section~\(f\) is continuous if and only if, for any
  \(\chi\in\widehat{G}\) and \(\varepsilon>0\), we can find a
  neighborhood~\(U\) of~\(\chi\) and \(a\in\W\) with
  \(\norm{\hat{a}(\omega)-f(\omega)}<\varepsilon\) for all \(\omega\in
  U\).  For a continuous section~\(f\) with compact support, we can find
  a finite covering by such open subsets \(U_1,\dotsc,U_n\) and
  \(a_1,\dotsc,a_n\in\W\) with
  \(\norm{\hat{a}_j(\omega)-f(\omega)}<\varepsilon\) for \(\omega\in
  U_j\).  In addition, for any fixed neighborhood~\(V\) of \(\supp f\),
  we can assume \(U_1,\dotsc,U_n\subseteq V\).  We can find continuous
  functions \(\phi_1,\dotsc,\phi_n\colon \widehat{G}\to[0,1]\) with
  \(\supp \phi_j\subseteq U_j\) and \(\phi_1+\dotsb+\phi_n=1\) on
  \(\supp f\).  We get
  \[
  \left\lVert f(\omega)- \sum_{i=1}^n \phi_i(\omega)\hat{a}_i(\omega)
  \right\rVert <\varepsilon
  \]
  for all \(\omega\in \bigcup U_j \supseteq \supp f\), and the sum
  vanishes outside \(\bigcup U_j\).  This shows that~\(f\) belongs to
  the closed \(\CONT_c(\widehat{G})\)-submodule generated by \(\{\hat{a}
  \mid a\in\W\}\).  The last assertion now follows because the Fourier
  transform maps \(\IG\) to a dense subspace of \(\CONT_c(\widehat{G})\)
  by Lemma \ref{lem:IG} and the multiplication maps
  \(\CONT_c(\widehat{G})\to \CONT_c(\Bundle{B})\), \(g\mapsto g\cdot
  \hat{a}\) are continuous for all \(a\in\W\).%\qed
\end{proof}

\begin{lem}
  \label{lem:iso_Fell}
  The isomorphism \(C^*\bigl(\Bundle{B}(\W)\bigr) \to A\) constructed in
  \textsection11 of~\cite{exel1} maps an \(\Ls^1\)-section
  \((b_\chi)_{\chi\in\widehat{G}}\) to \(\int^\su_{\widehat{G}} b_\chi
  \dd{\chi}\).
\end{lem}

\begin{proof}
  Sections of the form \(\widehat{g*a} = \hat{g}\cdot\hat{a}\) with
  \(a\in\W\) and \(g\in\IG\) span a dense subspace of
  \(\CONT_c(\Bundle{B})\) and hence of \(\Ls^1(\Bundle{B})\) by Lemma
  \ref{lem:Fell_cont_section} (compare Proposition 11.9
  in~\cite{exel1}).  Hence it suffices to check the assertion for
  sections of this form.  The Fourier inversion formula (Proposition
  \ref{pro:Fourier-inversion}) yields
  \[
  \int^\su_{\widehat{G}} \widehat{g*a}(\chi) \dd{\chi} = g*a.
  \]
  The proof of Proposition 11.10 in~\cite{exel1} shows that this is the
  image of \(\widehat{g*a}\) under the isomorphism
  \(C^*\bigl(\Bundle{B}(\W)\bigr)\to A\).  (The only point that makes
  this complicated is that Exel represents the crossed product
  \(\Cross{G}{A}\) on a Hilbert space instead of the Hilbert module
  \(L^2(G,A)\).)%\qed
\end{proof}

Now we use the equivalence between relatively continuous subsets of
integrable and square-integrable elements to rephrase Exel's results in
terms of continuously square-integrable Hilbert modules.

Let~\(\E\) be a Hilbert \(B,G\)\brd{}module and let
\(\Relc\subseteq\E_\si\) be relatively continuous, complete, and dense
in~\(\E\).  We get a Fell bundle
\[
\Bundle{B}(\E,\Relc) \defeq \Bundle{B}(\Comp(\E),\W_\Relc)
\]
with \(\W_\Relc\subseteq\Comp(\E)_\ii\) as in \eqref{eq:W_Relc}.  Using
\eqref{eq:E_x_ketbra}, we identify its fibers with
\begin{equation}
  \label{eq:Fell_bundle_fiber}
  \begin{aligned}[t]
    \Bundle{B}_\chi(\E,\Relc)
    &= \cspn {}\{E_\chi(\ket{\xi}\bra{\eta})\mid
    \xi,\eta\in\Relc\}\\
    &= \cspn {}\{\KET{\xi}\circ M_\chi\circ \BRA{\eta} \mid
    \xi,\eta\in\Relc\}\\
    &= \cspn {}\{\zeta\circ M_\chi\circ \upsilon^* \mid
    \zeta,\upsilon\in\F(\E,\Relc)\},
  \end{aligned}
\end{equation}
where~\(\F(\E,\Relc)\) is the Hilbert \(\Cross{G}{B}\)-module associated
to \((\E,\Relc)\) as in \eqref{eq:F_E_R}.  The algebraic structure of
the Fell bundle is inherited from \(\Bound(\E) =
\Mult\bigl(\Comp(\E)\bigr)\) using \(\Bundle{B}_\chi(\E,\Relc)\subseteq
\Bound(\E)\); the spaces of continuous sections are generated by
sections of the form \(\chi\mapsto \zeta\circ M_\chi\circ \upsilon^*\)
for \(\zeta,\upsilon\in\F(\E,\Relc)\) (compare Lemma
\ref{lem:Fell_cont_section}).

\begin{cor}
  \label{cor:fixed_as_fiber}
  The generalized fixed-point algebra \(\Fix(\E,\Relc)\) is the fiber at
  the unit element in the Fell bundle \(\Bundle{B}(\E,\Relc)\).
\end{cor}

\begin{rem}
  \label{rem:direct_Fell_from_Hmodule}
  Given any concrete Hilbert \(\Cross{G}{B}\)-module
  \[
  \F\subseteq \Bound^G(L^2(G,B),\E)
  \]
  (see~\cite{ralf1}), it can be checked directly that the last
  description in \eqref{eq:Fell_bundle_fiber} yields a Fell
  bundle~\(\Bundle{B}\) over~\(\widehat{G}\) and a canonical
  \(G\)\nbd{}equivariant isomorphism \(C^*(\Bundle{B}) \cong
  \Comp(\E)\).  We omit this argument because it is quicker (but less
  beautiful) to cite Exel's work.
\end{rem}

\begin{lem}
  \label{lem:dense_Relc_Fell_fiber}
  Let \((\E,\Relc)\) be a continuously square-integrable Hilbert
  \(B,G\)\brd{}module and let \(\Relc_0\subseteq\Relc\) be dense in the
  \(\si\)-norm.  Then
  \[
  \Bundle{B}_\chi(\E,\Relc)
  = \cspn {}\{\KET{\xi} \circ M_\chi \circ \BRA{\eta} \mid
  \xi,\eta\in\Relc_0\}
  \qquad\text{for all \(\chi\in\widehat{G}\).}
  \]
\end{lem}

\begin{proof}
  Since the maps \(\xi\mapsto \KET{\xi}, \BRA{\xi}\) are continuous for
  the \(\si\)-norm and the operator norm, this follows from
  \eqref{eq:Fell_bundle_fiber}.%\qed
\end{proof}

If \(\E=A\) for a \(G\)\nbd{}\Cstar{}algebra~\(A\), our construction
yields a continuous spectral decomposition of~\(A\) for any complete,
relatively continuous subset \(\Relc\subseteq A_\si\) that is dense
in~\(A\).

\begin{lem}
  \label{lem:Fell_from_Relc_functorial}
  The construction \((A,\Relc)\mapsto \Bundle{B}(A,\Relc)\) is a functor
  from the category of continuously square-integrable
  \(G\)\nbd{}\Cstar{}algebras to the category of Fell bundles
  over~\(\widehat{G}\).  In particular, if \((A_1,\Relc_1)\) and
  \((A_2,\Relc_2)\) are isomorphic, so are \(\Bundle{B}(A_1,\Relc_1)\)
  and \(\Bundle{B}(A_2,\Relc_2)\).
\end{lem}

We recall during the proof what the morphisms in both categories are.
This also specifies the \emph{iso}morphisms.

\begin{proof}
  Let \((A_t,\Relc_t)\) for \(t=0,1\) be two continuously
  square-integrable \(G\)\nbd{}\Cstar{}algebras.  Let
  \(\Bundle{B_t}\defeq \Bundle{B}(A_t,\Relc_t)\) be the associated Fell
  bundles.  A morphism \((A_1,\Relc_1)\to (A_2,\Relc_2)\) is a
  \(G\)\nbd{}equivariant \Star{}homomorphism \(f\colon A_1\to A_2\) with
  \(f(\Relc_1)\subseteq\Relc_2\).  We want to use such a morphism to
  construct a morphism of Fell bundles \(\Bundle{B_1}\to \Bundle{B_2}\).
  This is a continuous map \(\phi\colon \Bundle{B_1}\to\Bundle{B_2}\)
  that preserves the norm, product and involution, and such that
  \(\phi|_{\Bundle{B}_\chi}\) is a linear map from
  \(\Bundle{B}_{1,\chi}\) to \(\Bundle{B}_{2,\chi}\) for all
  \(\chi\in\widehat{G}\).  We can describe it equivalently as a family
  of linear maps \(\phi_\chi\colon \Bundle{B}_{1,\chi}\to
  \Bundle{B}_{2,\chi}\) for all \(\chi\in\widehat{G}\), with the
  following properties:
  \begin{itemize}
  \item \(\phi\) is compatible with the multiplication:
    \(\phi_{\chi_1\chi_2}(b_1 \cdot b_2) = \phi_{\chi_1}(b_1)\cdot
    \phi_{\chi_2}(b_2)\) for all \(\chi_1,\chi_2\in\widehat{G}\),
    \(b_1\in \Bundle{B}_{1,\chi_1}\), \(b_2\in \Bundle{B}_{1,\chi_2}\);

  \item \(\phi\) is compatible with the involution:
    \(\phi_{\cnj{\chi}}(b^*) = \phi_\chi(b)^*\) for all
    \(\chi\in\widehat{G}\), \(b\in \Bundle{B}_{1,\chi}\);

  \item \(\phi\) is continuous:
    \(\bigl(\phi_\chi(b_\chi)\bigr)_{\chi\in\widehat{G}}\) is a
    continuous section of~\(\Bundle{B_2}\) for any continuous section
    \((b_\chi)_{\chi\in\widehat{G}}\) of~\(\Bundle{B_1}\); it suffices
    to require this for a generating set of continuous sections (compare
    Lemma \ref{lem:Fell_cont_section}).

  \end{itemize}
  The first two conditions above imply \(\norm{\phi_\chi}\le1\) for all
  \(\chi\in\widehat{G}\).

  We do not require~\(f\) to be essential, that is, we allow
  \(f(A_1)\cdot A_2\neq A_2\).  Thus~\(f\) need not extend to a strictly
  continuous \Star{}homomorphism between the multiplier algebras.  To
  circumvent this issue, we work with the bidual von Neumann algebras;
  it is clear that~\(f\) induces a weakly continuous \Star{}homomorphism
  \(f''\colon A_1''\to A_2''\).  We have \(\Mult(A)\subseteq A''\) for
  any \Cstar{}algebra~\(A\), and the weak topology on~\(A''\) is weaker
  than the strict topology on \(\Mult(A)\).  Therefore, if an
  unconditional integral exists in the strict topology on \(\Mult(A)\),
  it also exists in the weak topology on~\(A''\).

  The map \(f\colon A_1\to A_2\) induces maps between the homogeneous
  subspaces \(\tilde{f}_\chi\colon (A_1'')_\chi \to (A_2'')_\chi\) ---
  defined as in \eqref{eq:chi_homogeneous} --- which automatically
  satisfy the first two algebraic conditions for a morphism of Fell
  bundles.  Since \(f(\Relc_1)\subseteq \Relc_2\) and since~\(f''\)
  preserves weak unconditional integrals, \(\tilde{f}_\chi\) maps
  \(E_\chi(a\cdot b^*)\) to \(E_\chi\bigl(f(a)\cdot f(b)^*\bigr)\) for
  all \(a,b\in\Relc_1\).  Hence we get maps \(\Bundle{B}_{1,\chi} \to
  \Bundle{B}_{2,\chi}\), which form a morphism of Fell bundles; the
  algebraic conditions are clear, the continuity condition holds because
  the generating continuous sections \(\chi\mapsto E_\chi(a\cdot b^*)\)
  for~\(\Bundle{B_1}\) are mapped to the continuous sections
  \(\chi\mapsto E_\chi(f(a)\cdot f(b)^*)\) of~\(\Bundle{B_2}\).%\qed
\end{proof}

\subsection{The inverse construction}
\label{sec:from_Fell_to_Relc}

Now we extract a relatively continuous \(G\)\nbd{}\Cstar{}algebra
\((A,\Relc)\) from a Fell bundle~\(\Bundle{B}\) over~\(\widehat{G}\).
Since this construction should be inverse to the one in
\textsection\ref{sec:construct_Fell}, we let \(A\defeq C^*(\Bundle{B})\)
with the dual action of~\(G\).

We view the fibers~\(\Bundle{B}_\chi\) of~\(\Bundle{B}\) as subspaces of
the multiplier algebra \(\Mult(A)\) via
\begin{equation}
  \label{eq:Fell_fiber_mult}
  (b_\chi\cdot f)(\eta) = b_\chi\cdot f(\chi^{-1}\eta),
  \qquad
  (f\cdot b_\chi)(\eta) = f(\eta\chi^{-1})\cdot b_\chi
\end{equation}
for \(f\in\Ls^1(\Bundle{B})\); this yields isometric embeddings
\(\Bundle{B}_\chi \to \Mult_\chi\bigl(C^*(\Bundle{B})\bigr)\) by
\cite[VIII.5.8]{fell_doran}.

\begin{lem}
  \label{lem:CONT_c_si}
  Let \(\Bundle{B} = (\Bundle{B}_\chi)_{\chi\in\widehat{G}}\) be a Fell
  bundle over~\(\widehat{G}\).  Then \(\CONT_c(\Bundle{B})\) consists of
  square-integrable elements of \(C^*(\Bundle{B})\).  Furthermore, we
  have
  \begin{equation}
    \label{eq:E_chi_for_CONT_c}
    E_\chi(ab^*) = (ab^*)(\chi)
  \end{equation}
  for all \(a,b\in\CONT_c(\Bundle{B})\), \(\chi\in\widehat{G}\), and
  \begin{equation}
    \label{eq:si-norm_Fell}
    \norm{f}_\si
    = \norm{f}_{C^*(\Bundle{B})}
    + \left\lVert \int_G f(\chi)\cdot f(\chi)^* \dd{\chi}
    \right\rVert_{\Bundle{B_1}}^{\nicefrac{1}{2}}
  \end{equation}
  for all \(f\in\CONT_c(\Bundle{B})\).
\end{lem}

\begin{proof}
  Theorem 5.5 in~\cite{exel2} asserts that elements of the form \(f^*f\)
  with \(f\in\CONT_c(\Bundle{B})\) are integrable with respect to the
  dual action and satisfy \(E_\chi(f^*f) = (f^*f)(\chi)\) for all
  \(\chi\in\widehat{G}\).  Since \(\CONT_c(\Bundle{B})\) is invariant
  under the involution, we may replace \(f^*f\) by \(ff^*\) here; we get
  \eqref{eq:E_chi_for_CONT_c} by polarization.  The integrability of
  \(\CONT_c(\Bundle{B})^2\) is equivalent to \(\CONT_c(\Bundle{B})
  \subseteq C^*(\Bundle{B})_\si\) by Proposition
  \ref{pro:Comp_integrable}.  Notice that \(ff^*=\ket{f}\bra{f}\) if we
  identify~\(A\) with the \Cstar{}algebra of compact operators on~\(A\)
  viewed as a Hilbert module over itself.

  Finally, the following computation implies \eqref{eq:si-norm_Fell}:
  \begin{multline*}
    \norm{\KET{f}}^2
    = \norm{\KET{f}\BRA{f}}_{\Mult(C^*(\Bundle{B}))}
    = \norm{E_1(f*f^*)}_{\Mult(C^*(\Bundle{B}))}\\
    = \norm{(f*f^*)(1)}_{\Mult(C^*(\Bundle{B}))}
    = \norm{(f*f^*)(1)}_{\Bundle{B_1}}\\
    = \left\lVert \int_{\widehat{G}} f(\chi)\cdot f^*(\chi^{-1}) \dd{\chi}
    \right\rVert_{\Bundle{B_1}}
    = \left\lVert \int_{\widehat{G}} f(\chi)\cdot f(\chi)^* \dd{\chi}
    \right\rVert_{\Bundle{B_1}}.
  \end{multline*}
  Here we use \eqref{eq:E_x_ketbra}, \eqref{eq:E_chi_for_CONT_c}, and
  that the embedding \(\Bundle{B_1}\to\Mult\bigl(C^*(B)\bigr)\) is
  isometric.%\qed
\end{proof}

\begin{prop}
  \label{pro:CONT_c_relc}
  Let \(\Bundle{B} = (\Bundle{B}_\chi)_{\chi\in\widehat{G}}\) be a Fell
  bundle over~\(\widehat{G}\).  Then the \(\si\)-norm closure
  \(\clos{\CONT_c(\Bundle{B})}^\si\) is a relatively continuous,
  complete, square-integrable, dense subspace of \(C^*(\Bundle{B})\), so
  that \(\bigl(C^*(\Bundle{B}), \clos{\CONT_c(\Bundle{B})}^\si\bigr)\)
  is a continuously square-integrable \(G\)\nbd{}\Cstar{}algebra.
\end{prop}

\begin{proof}
  Let \(\Relc\defeq \CONT_c(\Bundle{B})\) and let~\(\W_\Relc\) be the
  corresponding set of integrable elements as in \eqref{eq:W_Relc}.
  Proposition 10.2 in~\cite{exel1} asserts that~\(\W_\Relc\) is
  relatively continuous as a set of integrable elements.  By Proposition
  \ref{pro:compare_si_int}, this is equivalent to the relative
  continuity of~\(\Relc\) as a set of square-integrable elements.
  Since~\(\Relc\) is a dense, \(G\)\nbd{}invariant subalgebra of
  \(C^*(\Bundle{B})\), its completion agrees with its \(\si\)-norm
  closure by Corollary \ref{cor:completion}.%\qed
\end{proof}

\begin{defn}
  \label{def:csi_from_Fell}
  We define \(A(\Bundle{B})\defeq C^*(\Bundle{B})\) and
  \(\Relc(\Bundle{B}) \defeq \clos{\CONT_c(\Bundle{B})}^\si\) for a Fell
  bundle~\(\Bundle{B}\) over~\(\widehat{G}\).
\end{defn}

We get a functor \(\Bundle{B}\mapsto \bigl(A(\Bundle{B}),
\Relc(\Bundle{B})\bigr)\) from the category of Fell bundles to the
category of continuously square-integrable \(G\)\nbd{}\Cstar{}algebras.
Any morphism of Fell bundles \(\Bundle{B_1}\to\Bundle{B_2}\) induces a
\Star{}homomorphism \(C^*(\Bundle{B_1})\to C^*(\Bundle{B_2})\), which is
\(G\)\nbd{}equivariant for the dual actions and restricts to a map
\(\CONT_c(\Bundle{B_1})\to \CONT_c(\Bundle{B_2})\); it is contractive
for the \(\si\)-norm by \eqref{eq:si-norm_Fell}, so that it maps
\(\Relc(\Bundle{B_1}) = \clos{\CONT_c(\Bundle{B_1})}^\si\) to
\(\Relc(\Bundle{B_2}) = \clos{\CONT_c(\Bundle{B_2})}^\si\).

The Hilbert module \(\F \defeq
\F\bigl(A(\Bundle{B}),\Relc(\Bundle{B})\bigr)\) over
\(\Cross{G}{A(\Bundle{B})}\) associated to \(\bigl(A(\Bundle{B}),
\Relc(\Bundle{B})\bigr)\) already appears in \textsection10
of~\cite{exel1}, but in disguise.  It is the dual \(\mathcal{X}^*\defeq
\Comp(\mathcal{X},\Bundle{B_1})\) of the imprimitivity
bimodule~\(\mathcal{X}\) that implements a Morita--Rieffel equivalence
between \(\Bundle{B_1}\) and an ideal in \(\Cross{G}{C^*(\Bundle{B})}\).
Since the dualization exchanges left and right, the
\(\Bundle{B_1}\)\nbd{}valued inner product \(\int a(\chi)^* b(\chi)
\dd{\chi}\) in~\cite{exel1} becomes \(\int a(\chi) b(\chi)^*
\dd{\chi}\).

\subsection{The main theorem}
\label{sec:main}

We have constructed a functor from continuously square-integrable
\(G\)\nbd{}\Cstar{}algebras to Fell bundles over~\(\widehat{G}\) and a
functor in the opposite direction.

\begin{thm}
  \label{the:main}
  These two functors are inverse to each other (up to natural
  isomorphism).  That is, we have a natural isomorphism
  \[
  (A_0,\Relc_0) \cong \Bigl(A\bigl(\Bundle{B}(A_0,\Relc_0)\bigr),
  \Relc\bigl(\Bundle{B}(A_0,\Relc_0)\bigr)\Bigr)
  \]
  for any continuously square-integrable \(G\)\nbd{}\Cstar{}algebra
  \((A_0,\Relc_0)\), and a natural isomorphism
  \[
  \Bundle{B_0} \cong
  \Bundle{B}\bigl(A(\Bundle{B_0}),\Relc(\Bundle{B_0})\bigr)
  \]
  for any Fell bundle~\(\Bundle{B_0}\) over~\(\widehat{G}\).

  Hence the categories of continuously square-integrable
  \(G\)\nbd{}\Cstar{}algebras and of Fell bundles over~\(\widehat{G}\)
  are equivalent.
\end{thm}

Before we prove this, we use Theorem \ref{the:main} to classify
continuous spectral decompositions of a given
\(G\)\nbd{}\Cstar{}algebra~\(A\).  First we must discuss the
\emph{equivalence} of continuous spectral decompositions.

Recall that a continuous spectral decomposition is a Fell
bundle~\(\Bundle{B}\) with an isomorphism \(\phi\colon
C^*(\Bundle{B})\to A\).  The issue is whether~\(\phi\) is part of the
data or not.  If~\(\phi\) is part of the data, then Theorem
\ref{the:main} shows that isomorphism classes of continuous spectral
decompositions of~\(A\) correspond bijectively to relatively continuous,
square-integrable, complete dense subspaces \(\Relc\subseteq A\); that
is, two such subspaces yield isomorphic spectral decompositions if and
only if they are \emph{equal}.  The reason is that an isomorphism
\(\Bundle{B_1} \to \Bundle{B_2}\) that is compatible with the
isomorphisms to~\(A\) induces the identity map on \(A\cong
C^*(\Bundle{B_1}) \cong C^*(\Bundle{B_2})\).

Now let us disregard~\(\phi\) and consider two spectral decompositions
\(\phi_1\colon C^*(\Bundle{B})\to A\) and \(\phi_2\colon
C^*(\Bundle{B})\to A\) with the same bundle~\(\Bundle{B}\) to be
equivalent.  The automorphism \(\gamma\defeq \phi_2\circ \phi_1^{-1}\)
of~\(A\) maps \(\phi_1\bigl(\CONT_c(\Bundle{B})\bigr)\) to
\(\phi_2\bigl(\CONT_c(\Bundle{B})\bigr)\).  Therefore, their
\(\si\)-norm closures \(\Relc_1,\Relc_2\subseteq A\) are
\(\Aut(A)\)-conjugate in the following sense:

\begin{defn}
  \label{def:auto}
  Let \(\Aut(A)\) be the group of \(G\)\nbd{}equivariant \Cstar{}algebra
  automorphisms of~\(A\).  We call \(\Relc_1,\Relc_2\subseteq A\)
  \emph{\(\Aut(A)\)-conjugate} if there is \(\gamma\in\Aut(A)\) with
  \(\gamma(\Relc_1)=\Relc_2\).
\end{defn}

\begin{cor}
  \label{cor:main}
  Let~\(A\) be a \(G\)\nbd{}\Cstar{}algebra.  Isomorphism classes of
  Fell bundles~\(\Bundle{B}\) over~\(\widehat{G}\) for which there
  exists a \(G\)\nbd{}\Cstar{}algebra isomorphism \(C^*(\Bundle{B})
  \cong A\) correspond bijectively to \(\Aut(A)\)-conjugacy classes of
  relatively continuous, square-integrable, complete dense subspaces
  \(\Relc\subseteq A\).
\end{cor}

\subsection{Proof of the main theorem}
\label{sec:proof_main}

Let~\(\Bundle{B}\) be a Fell bundle over~\(\widehat{G}\).  We want to
construct a natural isomorphism \(\Bundle{B}\cong
\Bundle{B}(C^*(\Bundle{B}),\clos{\CONT_c(\Bundle{B})}^\si)\).

\begin{lem}
  \label{lem:square_dense}
  The subset \(\{\xi(\chi) \mid
  \xi\in\CONT_c(\Bundle{B})*\CONT_c(\Bundle{B})\}\) is dense
  in~\(\Bundle{B}_\chi\) for all \(\chi\in\widehat{G}\), and
  \(\CONT_c(\Bundle{B})*\CONT_c(\Bundle{B})\) is dense in
  \(\CONT_c(\Bundle{B})\).
\end{lem}

\begin{proof}
  It is clear that \(\xi(\chi)\) with \(\xi\in\CONT_c(\Bundle{B})\) are
  dense in~\(\Bundle{B}_\chi\) (see Remark II.13.19
  in~\cite{fell_doran}).  Remark VIII.5.12 in~\cite{fell_doran} asserts
  that \(\CONT_c(\Bundle{B})*\CONT_c(\Bundle{B})\) is dense in
  \(\CONT_c(\Bundle{B})\) for the inductive limit topology.%\qed
\end{proof}

We identify~\(\Bundle{B}_\chi\) with a subspace of
\(\Mult_\chi\bigl(C^*(\Bundle{B})\bigr)\) as in
\eqref{eq:Fell_fiber_mult}.

\begin{lem}
  \label{lem:004}
  We have \(\Bundle{B}_\chi = \clos{\{\hat{a}(\chi)\mid a\in
    \CONT_c(\Bundle{B})*\CONT_c(\Bundle{B})\}}\).
\end{lem}

\begin{proof}
  Recall that \(\CONT_c(\Bundle{B})\subseteq C^*(\Bundle{B})_\si\), so
  that \(\hat{a}(\chi) = E_\chi(a)\) is well-defined for \(a\in
  \CONT_c(\Bundle{B})*\CONT_c(\Bundle{B})\).  The assertion now follows
  from \eqref{eq:E_chi_for_CONT_c} and Lemma
  \ref{lem:square_dense}.%\qed
\end{proof}

\begin{prop}
  Let \(A \defeq C^*(\Bundle{B})\) and \(\Relc \defeq
  \clos{\CONT_c(\Bundle{B})}^\si\).  The inclusions
  \[
  \Bundle{B}_\chi\hookrightarrow \Mult_\chi(A) \hookleftarrow
  E_\chi(\W_\Relc)
  \]
  induce a natural isomorphism of Fell bundles
  \(\Bundle{B}\cong\Bundle{B}(A,\Relc)\).
\end{prop}

\begin{proof}
  In order to ensure that our construction is natural for all morphisms
  of Fell bundles, we replace \(\Mult(A)\) by the bidual \(A''\supseteq
  \Mult(A)\), which is functorial in complete generality.  We have
  described morphisms of Fell bundles in the proof of Lemma
  \ref{lem:Fell_from_Relc_functorial}.

  The inclusions into~\(A''_\chi\) identify
  \(\Bundle{B}_\chi\bigl(C^*(\Bundle{B}),
  \clos{\CONT_c(\Bundle{B})}^\si\bigr)\) and \(\Bundle{B}_\chi\) as
  Banach spaces by Lemmas \ref{lem:dense_Relc_Fell_fiber} and
  \ref{lem:004} and \eqref{eq:E_x_ketbra}.  The resulting family of
  isomorphisms \(\Bundle{B}_\chi \cong \Bundle{B}(A,\Relc)_\chi\) is
  compatible with multiplication and involution because both are
  inherited from \(\Mult(A)\subseteq A''\).  It remains to check that
  both Fell bundles carry the same topology.  By Lemma \ref{lem:004},
  sections of the form \(f*g^*\) with \(f,g\in\CONT_c(\Bundle{B})\)
  span a dense subspace in \(\CONT_c(\Bundle{B})\).  They span a dense
  subspace in \(\CONT_c\bigl(\Bundle{B}(A,\Relc)\bigr)\) as well by
  \eqref{eq:E_chi_for_CONT_c}.  Now we use that the topology of a Fell
  bundle is determined by its space of compactly supported continuous
  sections.%\qed
\end{proof}

This proves one half of our main theorem (Theorem \ref{the:main}).

Now let \((A,\Relc)\) be a continuously square-integrable
\(G\)\nbd{}\Cstar{}algebra.  We abbreviate \(\Bundle{B}\defeq
\Bundle{B}(A,\Relc)\).  We are going to show that
\[
(A,\Relc) \cong
\bigl(C^*(\Bundle{B}), \clos{\CONT_c(\Bundle{B})}^\si\bigr).
\]
An isomorphism \(\phi\colon C^*(\Bundle{B})\overset{\cong}\to A\) has
already been constructed by Ruy Exel; it is described in Lemma
\ref{lem:iso_Fell}.  As a \(G\)\nbd{}equivariant \Star{}isomorphism, it
must restrict to an isometric isomorphism \(C^*(\Bundle{B})_\si \cong
A_\si\).  It remains to check that~\(\phi\) maps
\(\clos{\CONT_c(\Bundle{B})}^\si\) onto~\(\Relc\).

Let \(\IG\) be as in Definition \ref{def:IG}, and define~\(\W_\Relc\) as
in \eqref{eq:W_Relc}.  The Fourier inversion theorem (Proposition
\ref{pro:Fourier-inversion}) and Lemma \ref{lem:iso_Fell} imply
\[
\phi(\widehat{g*a}) = \phi(\hat{g}\cdot \hat{a}) = g*a
\qquad\text{for all \(a\in\W_\Relc\), \(g\in\IG\).}
\]
Lemma \ref{lem:Fell_cont_section} asserts that the linear span of
sections of the form \(\widehat{g*a}\) is dense in
\(\CONT_c(\Bundle{B})\) in the inductive limit topology and hence in the
\(\si\)-norm.  Thus~\(\phi\) maps \(\clos{\CONT_c(\Bundle{B})}^\si\)
onto the \(\si\)-norm closed linear span of \(\IG*\W_\Relc\).

\begin{lem}
  \label{lem:W_versus_Relc}
  The \(\si\)-norm closed linear span of \(\IG*\W_\Relc\) is~\(\Relc\).
\end{lem}

\begin{proof}
  Since~\(\Relc\) is a right ideal, it contains \(\W_\Relc = \Relc\cdot
  \Relc^*\).  Since~\(\Relc\) is \(G\)\nbd{}invariant and \(\si\)-norm
  closed, it also contains the \(\si\)-norm closed linear span of
  \(\IG*\W_\Relc\).

  Since~\(G\) acts continuously on~\(\Relc\), \(\Relc\) is an essential
  \(L^1(G)\)-module, that is, \(L^1(G)*\Relc\) is dense in~\(\Relc\).
  So is \(\IG*\Relc\) because \(\IG\) is dense in \(L^1(G)\) by Lemma
  \ref{lem:IG}.  Thus it remains to prove that~\(\W_\Relc\) is dense
  in~\(\Relc\).  This follows because~\(\Relc\) is an essential right
  \(A\)\nbd{}module and~\(\Relc^*\) is dense in~\(A\), so that
  \(\W_\Relc = \Relc\cdot \Relc^*\) is dense in~\(\Relc\).%\qed
\end{proof}

Since~\(\phi\) is isometric in the \(\si\)-norm and~\(\Relc\) is
complete, Lemma \ref{lem:W_versus_Relc} implies that~\(\phi\) maps
\(\clos{\CONT_c(\Bundle{B})}^\si\) onto~\(\Relc\).  This finishes the
proof of our main theorem (Theorem \ref{the:main}).

\section{Examples}
\label{sec:examples}

\subsection{Spectral decomposition for algebras of compact operators}
\label{sec:ex_Hilbert}

Let~\(G\) be a second countable locally compact Abelian group and
let~\(\Hils\) be a separable \(G\)\nbd{}Hilbert space, that is, a
separable Hilbert space with a continuous unitary representation~\(\pi\)
of~\(G\).  We let~\(G\) act on \(A\defeq \Comp(\Hils)\) by conjugation:
\[
\alpha_g(T)\defeq \Ad(\pi_g)(T) \defeq \pi_g\circ T\circ \pi_g^{-1}.
\]
What are the continuous spectral decompositions of~\(A\)?

By our main theorem, these correspond to complete relatively continuous
subspaces \(\Relc\subseteq A_\si\) that are dense in~\(A\).  By Theorem
7.2 of~\cite{ralf1}, such subspaces of~\(A\) correspond bijectively to
complete relatively continuous subspaces \(\Relc\subseteq \Hils_\si\)
that are dense in~\(\Hils\).  The latter are already classified in
\textsection8 of~\cite{ralf1}.  We briefly recall this and apply it to
continuous spectral decompositions.

\begin{rem}
  \label{rem:proj_rep}
  Any continuous action of~\(G\) on \(\Comp(\Hils)\) comes from a
  \emph{projective} representation of~\(G\) on~\(\Hils\).  But we only
  consider honest Hilbert space representations here for simplicity.
\end{rem}

If a subspace~\(\Relc\) as above exists, then~\(\Hils\) must be
square-integrable in the sense of Rieffel~\cite{rieffel1}.  The
\(G\)\nbd{}equivariant stabilization theorem (Theorem 8.5 in
\cite{ralf2}) shows that~\(\Hils\) is a direct summand in an infinite
direct sum of regular representations \(\bigoplus_{n\in\Nset} L^2(G)\).
Equivalently, \(\pi\) integrates to a normal \Star{}representation of
the group von Neumann algebra of~\(G\).  The Fourier transform
identifies the latter with \(L^\infty(\widehat{G})\).

Any normal \Star{}representation of \(L^\infty(\widehat{G})\) is
unitarily equivalent to the representation by pointwise multiplication
operators on the Hilbert space of \(L^2\)\nbd{}sections of a
\emph{measurable field} of Hilbert spaces over~\(\widehat{G}\).
Measurable fields of Hilbert spaces
\((\Hils_\chi)_{\chi\in\widehat{G}}\) over~\(\widehat{G}\) are, in turn,
classified by their \emph{spectral multiplicity function}
\[
d\colon \widehat{G}\to \clos{\Nset} = \Nset\cup\{\infty\},
\qquad d(\chi)\defeq \dim\Hils_\chi.
\]
More precisely, two normal \Star{}representations of
\(L^\infty(\widehat{G})\) are unitarily equivalent if and only if their
spectral multiplicity functions agree almost everywhere.  Moreover, any
measurable function \(d\colon \widehat{G}\to\clos{\Nset}\) occurs as a
spectral multiplicity function.

Summing up, we get a bijection between unitary equivalence classes of
square-integrable representations of~\(G\) and equivalence classes of
measurable functions \(\widehat{G}\to\clos{\Nset}\) up to equality
almost everywhere \cite{ralf1}.

Isomorphism classes of \emph{continuously} square-integrable Hilbert
space representations \((\Hils,\Relc)\) correspond bijectively to
isomorphism classes of Hilbert modules over \(\Cross{G}{\Cset} \cong
\CONT_0(\widehat{G})\).  These Hilbert modules correspond to
\emph{continuous} fields of Hilbert spaces over~\(\widehat{G}\).

If we fix~\(\Hils\), then we get a bijection between complete relatively
continuous subspaces \(\Relc\subseteq\Hils_\si\) dense in~\(\Hils\) and
isomorphism classes of pairs \((\F,\phi)\), where~\(\F\) is a Hilbert
\(\CONT_0(\widehat{G})\)-module and~\(\phi\) is a \(G\)\nbd{}equivariant
unitary operator
\[
\phi\colon \F\otimes_{\CONT_0(\widehat{G})} L^2(G) \to \Hils.
\]
The equivalence class of the pair \((\F,\phi)\) corresponds to the
completion of the subspace \(\phi\bigl(\F\otimes
\CONT_c(G)\bigr)\subseteq\Hils_\si\).

In order to translate everything to measurable and continuous fields of
Hilbert spaces, we replace the left regular representation of~\(G\) on
\(L^2(G)\) by the unitarily equivalent representation on
\(L^2(\widehat{G})\) by pointwise multiplication operators:
\[
\gamma_t(\xi)(\chi)\defeq \cnj{\braket{\chi}{t}} \cdot \xi(\chi)
\qquad
\text{for \(\xi\in L^2(\widehat{G})\), \(\chi\in\widehat{G}\) and \(t\in G\).}
\]
This representation is equivalent to the left regular representation on
\(L^2(G)\) via the Fourier transform
\[
U\colon L^2(G)\to L^2(\widehat{G}),
\qquad
Uf(\chi) =
\hat{f}(\chi) \defeq \int_G f(t)\cdot\cnj{\braket{\chi}{t}}\dd{t},
\]
as defined in \eqref{eq:Fourier_scalar}.  The integrated form
of~\(\gamma_t\) is the \Star{}representation of \(\CONT_0(\widehat{G})
\cong C^*(G)\) by pointwise multiplication operators.  The functor
\[
\F\mapsto \F\otimes_{C^*(G)} L^2(G) \cong
\F\otimes_{\CONT_0(\widehat{G})} L^2(\widehat{G})
\]
is the forgetful functor that views a continuous field of Hilbert spaces
over~\(\widehat{G}\) as a measurable field.  The subspace
\(\Relc\subseteq\Hils \cong
L^2\bigl((\Hils_\chi)_{\chi\in\widehat{G}}\bigr)\) is the subspace of
continuous \(L^2\)-sections in~\(\Hils\) vanishing at infinity.

\begin{thm}
  \label{the:decompose_compacts}
  Let~\(\Hils\) be a square-integrable representation.  There is a
  bijection between isomorphism classes of continuous spectral
  decompositions of \(\Comp(\Hils)\) and continuous structures on the
  measurable field of Hilbert spaces underlying~\(\Hils\).

  Let \(d\colon\widehat{G}\to\clos{\Nset}\) be the spectral multiplicity
  function.  A continuous spectral decomposition exists if and only if
  there is a decreasing sequence of open subsets \((U_n)_{n\in\Nset}\)
  of~\(\widehat{G}\) with \(U_0=\widehat{G}\), \(d(\chi)=n\) for almost
  all \(\chi\in U_n\setminus U_{n+1}\), and \(d(\chi)=\infty\) for
  almost all \(\chi\in\bigcap_{n\in\Nset} U_n\).
\end{thm}

\begin{proof}
  By our Main Theorem (Theorem \ref{the:main}) and Theorem 7.2
  in~\cite{ralf1}, the continuous spectral decompositions of
  \(\Comp(\Hils)\) correspond to complete relatively continuous
  subspaces \(\Relc\subseteq\Hils_\si\) dense in~\(\Hils\).  The
  discussion above shows that these correspond to continuous structures
  on the measurable field underlying~\(\Hils\).  This yields the first
  statement.  Since measurable fields are classified by their dimension
  function, the second statement follows from the first one and the
  observation that dimension functions of \emph{continuous} fields of
  Hilbert spaces are lower semi-continuous.%\qed
\end{proof}

\begin{exmp}
  \label{exa:si_not_csi}
  Let \(S\subseteq\widehat{G}\) be a closed subset with non-zero measure
  and empty interior.  Then the subspace \(L^2(S)\) of functions in
  \(L^2(\widehat{G})\) vanishing outside~\(S\) is a non-zero
  \(G\)\nbd{}invariant Hilbert subspace.  Its spectral multiplicity
  function is the characteristic function of~\(S\).  It is easy to see
  that there is no lower semi-continuous function that is almost
  everywhere equal to it (see \textsection8 in \cite{ralf1}).  Hence
  Theorem \ref{the:decompose_compacts} shows that
  \(\Comp\bigl(L^2(S)\bigr)\) is a square-integrable
  \(G\)\nbd{}\Cstar{}algebra without continuous spectral decomposition.
  This answers negatively Questions 9.4 and 11.16 in~\cite{exel1}.
\end{exmp}

If \(A\defeq \Comp(\Hils)\) has a continuous spectral decomposition, it
is never unique because there are always several continuous structures
on a measurable field of Hilbert spaces.  Even if we work up to
\(\Aut(A)\)-conjugacy (see Definition \ref{def:auto}), the spectral
decomposition rarely becomes unique.  We only study the case of
\(\Hils=L^2(\widehat{G})\) for the sake of concreteness.

\begin{thm}
  \label{the:csi_on_hatG}
  The isomorphism classes of continuous spectral decompositions of\/
  \(\Comp\bigl(L^2(\widehat{G})\bigr)\) correspond bijectively to
  isomorphism classes of triples \((S,V,\psi)\), where
  \begin{itemize}
  \item \(S\) is an open subset of~\(\widehat{G}\) of full measure,

  \item \(V\) is a Hermitian complex line bundle on~\(S\), and

  \item \(\psi\) is a measurable section of~\(V\) with
    \(\abs{\psi(\chi)}=1\) for all \(\chi\in S\).

  \end{itemize}
  Two triples \((S_1,V_1,\psi_1)\) and \((S_2,V_2,\psi_2)\) are
  isomorphic if and only if \(S_1=S_2\) and there is an isomorphism
  \(\phi\colon V_1\cong V_2\) of Hermitian complex line bundles with
  \(\phi_*(\psi_1) = \psi_2\), that is,
  \(\phi\circ\psi_1(\chi)=\psi_2(\chi)\) for all \(\chi\in\widehat{G}\).

  Conjugacy classes of continuous spectral decompositions of
  \(\Comp\bigl(L^2(\widehat{G})\bigr)\) correspond bijectively to
  conjugacy classes of pairs \((S,V)\) with \(S\) and~\(V\) as above,
  where \((S_1,V_1)\) and \((S_2,V_2)\) are \emph{conjugate} if and only
  if there is some \(\chi\in\widehat{G}\) with \(\chi\cdot S_1=S_2\) and
  \(\chi^*(V_2) \cong V_1\), where~\(\chi^*\) means that we pull back
  the line bundle \(V_2\) on~\(S_2\) along the map \(\omega\mapsto
  \chi\cdot\omega\) to a line bundle on~\(S_1\).
\end{thm}

\begin{proof}
  Theorem \ref{the:decompose_compacts} reduces the problem to that of
  classifying continuous structures on the measurable field of Hilbert
  spaces over~\(\widehat{G}\) underlying \(L^2(\widehat{G})\).  Since
  its spectral multiplicity function is constant equal to~\(1\), we must
  consider continuous fields of Hilbert spaces
  \((\Hils_\chi)_{\chi\in\widehat{G}}\) with \(\dim \Hils_\chi = 1\) for
  almost all \(\chi\in\widehat{G}\).  Since \(\dim\Hils_\chi\) is lower
  semi-continuous, this implies \(\dim \Hils_\chi = 1_S\) for some
  \emph{open} subset \(S\subseteq\widehat{G}\) of \emph{full measure},
  that is, \(\abs{\widehat{G}\setminus S}=0\).

  If \(\omega\in S\), there is a non-zero continuous section of
  \((\Hils_\chi)_{\chi\in\widehat{G}}\) in a neighborhood of~\(\omega\).
  This provides a local trivialization of \((\Hils_\chi)\)
  near~\(\omega\) because \(\dim \Hils_\chi=1\) on~\(S\).  Thus our
  continuous field of Hilbert spaces is locally trivial on~\(S\).
  Equivalently, it is a Hermitian complex line bundle, that is, a
  complex line bundle with a continously varying family of inner
  products on the fibers.

  Let~\(V\) be a Hermitian complex line bundle over~\(S\).  Local
  trivializations of~\(V\) allow us to construct a \emph{measurable}
  section~\(\psi\) of~\(V\) with \(\abs{\psi}=1\) almost everywhere.
  Then \(f\mapsto f\cdot\psi\) defines a
  \(\widehat{G}\)\nbd{}equivariant unitary operator from
  \(L^2(\widehat{G})\) to the space \(L^2(V)\) of \(L^2\)\nbd{}sections
  of~\(V\) or, equivalently, an isomorphism between the underlying
  measurable fields of Hilbert spaces over~\(\widehat{G}\).  Conversely,
  any \(\widehat{G}\)\nbd{}equivariant unitary operator
  \(L^2(\widehat{G})\to L^2(V)\) is of this form for some~\(\psi\) as
  above.

  As a result, a triple \((S,V,\psi)\) as above specifies a continuous
  structure on the \(1\)\nbd{}dimensional constant measurable field of
  Hilbert spaces over~\(\widehat{G}\).  If two such triples yield the
  same structure, then \(S_1=S_2\) and \(V_1\cong V_2\) because the
  corresponding Hilbert modules over \(\CONT_0(\widehat{G})\) are
  isomorphic.  The isomorphism \(\phi\colon V_1\to V_2\) must satisfy
  \(\phi_*(\psi_1)= \psi_2\) in order to be compatible with the
  isomorphisms \(L^2(V_1)\cong L^2(\widehat{G}) \cong L^2(V_2)\).
  Conversely, such an isomorphism clearly identifies the resulting
  continuous structures on the measurable field underlying
  \(L^2(\widehat{G})\).

  Finally, we turn to the classification up to conjugacy.  It is
  well-known that any automorphism of \(\Comp(\Hils)\) has the form
  \(\Ad_u(T)\defeq uTu^*\) for some unitary operator
  \(u\in\Bound(\Hils)\), which is determined uniquely up to scalar
  multiples.  The automorphism \(\Ad_u\) is \(G\)\nbd{}equivariant if
  and only if \(\pi_t u\pi_t^{-1} = \chi(t)\cdot u\) for all \(t\in G\),
  for some function~\(\chi\) from~\(G\) to the unit circle in~\(\Cset\).
  Of course, \(\chi\) must be a character, so that we get
  \(u\in\Bound_\chi(\Hils)\) for some \(\chi\in\widehat{G}\).  These
  subspaces are easy to describe for \(\Hils=L^2(\widehat{G})\):

  \begin{lem}
    \label{lem:chi_homogeneous_Ltwo}
    Let \(T_\chi f(\omega) \defeq f(\chi^{-1}\cdot\omega)\), then
    \begin{equation}
      \label{eq:Bound_chi_concrete}
      \Bound_{\cnj{\chi}}(L^2\widehat{G})
      = L^\infty(\widehat{G})\cdot T_\chi
      \cong L^\infty(\widehat{G})
    \end{equation}
    where we represent \(L^\infty(\widehat{G})\) isometrically by
    pointwise multiplication operators.
  \end{lem}

  \begin{proof}
    It is easy to see that \(\gamma_t\circ T_\chi\circ\gamma_t^{-1} =
    \cnj{\braket{\chi}{t}} \cdot T_\chi\) for all
    \(\chi\in\widehat{G}\), \(t\in G\).  Hence
    \(T_\chi\in\Bound_{\cnj{\chi}}(L^2\widehat{G})\).  Since~\(T_\chi\)
    is unitary, we have
    \[
    u\in\Bound_{\cnj{\chi}}(L^2\widehat{G}) \quad\iff\quad
    u\cdot  T_\chi^* \in
    \Bound_1(L^2\widehat{G}) = \Bound^G(L^2\widehat{G});
    \]
    it is well-known that the latter agrees with
    \(L^\infty(\widehat{G})\).%\qed
  \end{proof}

  Since~\(T_\chi\) is unitary, the \(G\)\nbd{}equivariant automorphisms
  of \(\Comp(\Hils)\) are precisely the maps \(\Ad_{M_f\circ T_\chi}\)
  for some \(\chi\in\widehat{G}\) and some unitary measurable function
  \(f\in L^\infty(\widehat{G})\).

  The corresponding action on continuous structures on~\(\Hils\) is the
  obvious one: apply \(M_f\,T_\chi\) to the space of continuous
  sections.  Application of~\(M_f\) does not change \(S\) and~\(V\) and
  replaces~\(\psi\) by \(\psi\cdot f\).  Since any two measurable
  sections of norm~\(1\) are related in this fashion, the conjugacy
  class of a continuous structure is independent of~\(\psi\), so that we
  may drop this component.  Application of~\(T_\chi\) shifts functions
  by~\(\chi\), so that we replace \((S,\chi^*(V))\) by \((\chi\cdot
  S,V)\).%\qed
\end{proof}

Recall that the first Chern class classifies isomorphism classes of
complex line bundles on~\(S\) by the cohomology group \(H^2(S;\Zset)\).
The Hermitian inner product on a complex line bundle is unique up to
isomorphism.  Thus we can also rewrite our classification in terms of
pairs \((S,x)\) with~\(S\) as above and a cohomology class \(x\in
H^2(S;\Zset)\).

\begin{exmp}
  \label{exa:Ltwo_hatG}
  Consider the case \(G=\Zset\) with dual
  \(\widehat{G}=\Torus=\Rset/\Zset\).  Since~\(\Torus\) is
  \(1\)\nbd{}dimensional, subsets of~\(\Torus\) carry no non-trivial
  line bundles.  Hence continuous spectral decompositions of
  \(\Comp(L^2\Torus)\) correspond to pairs \((S,\psi)\) with an open
  subset \(S\subseteq \Torus\) of full measure and a unitary element
  \(\psi\in L^\infty(S)\).  Up to conjugacy, they are classified by open
  subsets \(S\subseteq\Torus\) of full measure, with \([S_1]=[S_2]\) if
  and only if there is \(\chi\in\Torus\) with \(S_1=\chi\cdot S_2\).
\end{exmp}

\begin{exmp}
  \label{exa:Ltwo_hatG_Zsquare}
  Consider the discrete group \(G=\Zset^2\) with dual
  \(\widehat{G}=\Torus^2\).  If \(S\subseteq\widehat{G}\) is a proper
  open subset, then \(H^2(S;\Zset)=0\) because~\(S\) is a
  \emph{non-compact} oriented \(2\)\nbd{}dimensional manifold.
  Hence~\(S\) supports no non-trivial line bundles in this case.  The
  conjugacy classes of continuous spectral decompositions of
  \(\Comp(L^2\Torus^2)\) of this kind are classified by proper open
  subsets of full measure up to translation as in Example
  \ref{exa:Ltwo_hatG}.

  But~\(\Torus^2\) carries non-trivial line bundles because
  \(H^2(\Torus^2;\Zset)=\Zset\).  We have \(\chi^*(V) \cong V\) for any
  line bundle~\(V\) on~\(\Torus^2\) and any \(\chi\in\Torus^2\)
  because~\(\Torus^2\) is path-connected.  Therefore, conjugate line
  bundles over~\(\Torus^2\) are isomorphic and non-isomorphic line
  bundles yield non-conjugate continuous spectral decompositions of
  \(\Comp(L^2\Torus^2)\).

  We can describe the non-trivial vector bundles on~\(\Torus^2\)
  explicitly.  We parametrize points in~\(\Torus^2\) by
  \(\Torus\times[0,1]\), with \((t,0)\sim(t,1)\) for all \(t\in\Torus\).
  A line bundle~\(V\) on~\(\Torus^2\) pulls back to a trivial bundle on
  \(\Torus\times[0,1]\) because the latter space is homotopy equivalent
  to the \(1\)\nbd{}dimensional space~\(\Torus\) and thus carries no
  non-trivial line bundles.  To reconstruct~\(V\) from a trivial line
  bundle on \(\Torus\times[0,1]\), we glue together the restrictions to
  trivial line bundles on \(\Torus\times\{0\}\) and
  \(\Torus\times\{1\}\).  This can be done by any continuous map
  \(\Torus\to \mathrm{U}(1) =\Torus\).  Since only the homotopy class of
  this map matters, it suffices to use the gluing functions \(z\mapsto
  z^n\).  Thus all complex line bundles on~\(\Torus^2\) are of the form
  \[
  V_n \defeq \Torus\times [0,1]\times \Cset /
  \text{\((z,1,x)\sim (z,0,z^n\cdot x)\) for all \(z\in\Torus\),
    \(x\in\Cset\)}
  \]
  for some \(n\in\Zset\).  This is the pedestrian way to identify the
  set of complex line bundles on~\(\Torus^2\) with \(\Zset\cong
  H^2(\Torus^2,\Zset)\).  The resulting spaces of continuous sections
  are
  \[
  \CONT(\Torus^2,V_n) =
  \{ f\in\CONT(\Torus\times[0,1]) \mid f(z,1) = z^n\cdot f(z,0) \}.
  \]
  We identify this with a subset of \(L^\infty(\widehat{G})\) using the
  unique measurable section~\(\psi\) with \(f(z,t)=1\) for
  \(z\in\Torus\), \(t\in [0,1)\).  This identifies
  \(\CONT(\Torus^2,V_n)\) with the space of all functions in
  \(L^\infty(\widehat{G})\) that are continuous as functions on
  \(\Torus\times [0,1)\) and \(\Torus\times (0,1]\) and satisfy
  \[
  \lim_{t\nearrow 1} f(z,t) = z^n \cdot \lim_{t\searrow 0} f(z,t).
  \]
\end{exmp}

Now we describe the continuous spectral decompositions of
\(\Comp(L^2\widehat{G})\) that result from a triple \((S,V,\psi)\) as in
Theorem \ref{the:csi_on_hatG}.

We can choose a \emph{continuous} unitary section~\(\psi\) of~\(V\) if
and only if~\(V\) is trivial.  If~\(\psi\) is continuous, then
\begin{equation}
  \label{eq:Relc_S}
  \Relc_S = L^2(\widehat{G})\cap \CONT_0(S).
\end{equation}
In general, \(\psi\) is discontinuous and we get
\begin{equation}
  \label{eq:Relc_S_V_psi}
  \Relc_{S,V,\psi} = \{f\in L^2(\widehat{G}) \mid
  \text{\(\psi\cdot f\) is a \(\CONT_0\)-section of~\(V\) on~\(S\)}\}.
\end{equation}

To describe the associated Fell bundle, we compute the operators
\(\KET{\xi}\).  We prefer to work with \(\KET{\xi}\circ U^*\colon
L^2(\widehat{G})\to L^2(\widehat{G})\), where
\[
U\colon L^2(G)\to L^2(\widehat{G}),
\qquad f\mapsto \hat{f}.
\]
is the Fourier transform.  Since \(\KET{\xi}\circ U^*\) is
\(G\)\nbd{}equivariant, it belongs to \(L^\infty(\widehat{G}) =
\Bound^G(L^2\widehat{G}) \subseteq \Bound(L^2\widehat{G})\).

If \(f\in\CONT_c(G)\) and \(\chi\in\widehat{G}\), then
\begin{multline*}
  (\KET{\xi}f)(\chi)
  = \int_G \gamma_t(\xi)(\chi) \cdot f(t) \dd{t}
  \\= \int_G \xi(\chi) \cdot f(t) \cdot  \cnj{\braket{\chi}{t}} \dd{t}
  = \xi(\chi) \cdot \hat{f}(\chi).
\end{multline*}
Thus \(\KET{\xi}(f) = \xi\cdot U(f)\) for all \(f\in\CONT_c(G)\).  This
extends to a bounded operator \(L^2(G)\to L^2(\widehat{G})\) if and only
if multiplication by~\(\xi\) is bounded on \(L^2(\widehat{G})\), if and
only if \(\xi\in L^\infty(\widehat{G})\).  Thus
\[
L^2(\widehat{G})_\si = L^\infty(\widehat{G})\cap L^2(\widehat{G}),
\]
and
\begin{equation}
  \label{eq:KET_xi_Fourierstar}
  \KET{\xi}\circ U^* = M_\xi
  \qquad \forall\; \xi\in L^\infty(\widehat{G})\cap L^2(\widehat{G}),
\end{equation}
where~\(M_\xi\) denotes the pointwise multiplication operator.
Moreover,
\begin{equation}
  \label{eq:ex_Hilbert_si_norm}
  \norm{\xi}_\si = \norm{\xi}_2 + \norm{\xi}_\infty.
\end{equation}

The isomorphism \(\Cross{G}{\Cset}\cong \CONT_0(\widehat{G})\) that we
have already used above is simply given by \(T\mapsto UTU^*\).  If
\(\xi,\eta\in L^\infty(\widehat{G})\cap L^2(\widehat{G})\), then
\[
\BRAKET{\xi}{\eta} = U^* M_\xi^* M_\eta U = U^* M_{\cnj{\xi}\eta} U.
\]
Thus \(\xi \rc \eta\) if and only if \(\cnj{\xi}\cdot \eta \in
\CONT_0(\widehat{G})\).  It is easy to verify with this criterion that
the subspaces \(\Relc_{S,V,\psi}\) in \eqref{eq:Relc_S_V_psi} are
relatively continuous, confirming our computations above.

For any relatively continuous subspace~\(\Relc\), the fibers
\(\Bundle{B}_\chi(A,\Relc)\) of the associated Fell bundle are
norm-closed linear subspaces of
\[
\Mult_\chi(A) = \Bound_\chi(L^2\widehat{G})
= M_{L^\infty(\widehat{G})} \circ T_\chi^*
\]
by Lemma \ref{lem:chi_homogeneous_Ltwo}.  Since the representation of
\(L^\infty(\widehat{G})\) by multiplication operators is isometric, we
may identify \(\Bundle{B}_\chi(A,\Relc)\) with a subspace of
\(L^\infty(\widehat{G})\).  We have \(T_\chi M_f = M_{\chi\bullet f}
T_\chi\) with \(\chi \bullet f(\omega)\defeq f(\chi^{-1}\omega)\) for
all \(\chi\in\widehat{G}\), \(f\in L^\infty(\widehat{G})\).  Hence
\begin{equation}
  \label{eq:multiply_Fell_concrete}
  (M_f T_\chi^*) \cdot (M_g T_\omega^*)
  = M_{f\cdot (\cnj{\chi} \bullet g)} T_{\chi\omega}^*,
  \qquad
  (M_f\,T_\chi^*)^* = M_{\chi\bullet \cnj{f}} T_{\cnj{\chi}}^*.
\end{equation}
This allows us to translate the multiplication and involution in
\(\Mult(A)\) to the corresponding subspaces of
\(L^\infty(\widehat{G})\).

The computation \(UM_\chi f(\omega) = f(\chi\omega) =
T_\chi^*Uf(\omega)\) for \(\chi,\omega\in\widehat{G}\) yields
\[
UM_\chi U^* = T_\chi^*\colon L^2(\widehat{G}) \to L^2(\widehat{G})
\]
for all \(\chi\in\widehat{G}\).  Combining this with
\eqref{eq:Fell_bundle_fiber}, \eqref{eq:KET_xi_Fourierstar}, and
\eqref{eq:multiply_Fell_concrete}, we get
\begin{align*}
  \Bundle{B}_\chi(A,\Relc) &=
  \cspn {}\{\KET{\xi}\circ U^*U M_\chi U^*U \circ\BRA{\eta} \mid
  \xi,\eta\in\Relc\}
  \\&= \cspn{} \{M_\xi T_\chi^* M_\eta^*\mid \xi,\eta\in \Relc\}
  \\&= \cspn{} \{M_{\xi\cdot(\cnj{\chi} \bullet \cnj{\eta})} T_\chi^*
  \mid \xi,\eta\in\Relc\}.
\end{align*}
Identifying \(\Bundle{B}_\chi(A,\Relc)\) with a subspace of
\(L^\infty(\widehat{G})\) as above, we get
\begin{equation}
  \label{eq:Fell_fiber_explicit}
  \Bundle{B}_\chi(A,\Relc) \cong
  \cspn{} \Relc\cdot (\cnj{\chi}\bullet\cnj{\Relc}).
\end{equation}

If~\(\Relc\) is given by \eqref{eq:Relc_S}, the right-hand side in
\eqref{eq:Fell_fiber_explicit} becomes
\[
\Bundle{B}_\chi(A,\Relc_S) \cong
\cspn{} \CONT_0(S)\cdot\CONT_0(\cnj{\chi} S)
= \CONT_0(S\cap \cnj{\chi} S).
\]
We may view sections of \(\Bundle{B} \defeq \Bundle{B}(A,\Relc_S)\) as
functions on \(\widehat{G}\times\widehat{G}\).  As one should expect,
continuity of \emph{sections} of this bundle simply corresponds to
continuity of the corresponding \emph{functions} on
\(\widehat{G}\times\widehat{G}\).  The algebraic operations can be
deduced from \eqref{eq:multiply_Fell_concrete}.

Next we consider the case \(\Relc_{S,V,\psi}\) where~\(V\) is the
trivial vector bundle, so that~\(\psi\) is a measurable scalar-valued
function.  Then \eqref{eq:Fell_fiber_explicit} yields
\[
\Bundle{B}_\chi(A,\Relc_{S,\psi}) \cong
\CONT_0(S\cap \cnj{\chi} S) \cdot
\cnj{\psi} \cdot(\cnj{\chi}\bullet\psi).
\]
Theorem \ref{the:csi_on_hatG} predicts that this continuous spectral
decomposition of~\(A\) is conjugate to the one for \(\psi=1\).
Concretely, this corresponds to the isomorphism of Fell bundles
\(\Bundle{B}(A,\Relc_{S,1}) \cong \Bundle{B}(A,\Relc_{S,\psi})\) given
by pointwise multiplication with the function \(\chi\mapsto \cnj{\psi}
\cdot(\cnj{\chi}\bullet\psi)\).

The case \(\Relc_{S,V,\psi}\) for a general Hermitian complex line
bundle~\(V\) is more complicated.  A useful general observation is that
the generalized fixed-point algebra is independent of \(V\)
and~\(\psi\):
\begin{equation}
  \label{eq:Fell_fiber_at_1}
  \Bundle{B_1}(A,\Relc_{S,V,\psi}) = \CONT_0(S).
\end{equation}
This is because the algebra bundle of endomorphisms of a complex line
bundle is always trivial, the identity section providing a nowhere
vanishing global section.  The generalized fixed-point algebras becomes
more complicated when we study \(\Comp(L^2\widehat{G}\oplus
L^2\widehat{G})\) (see \textsection8 in~\cite{ralf1}).

The other fibers of our Fell bundle can be described as follows:
\begin{multline}
  \label{eq:Fell_fiber_at_chi}
  \Bundle{B}_\chi(A,\Relc_{S,V,\psi})
  = \{f \in L^\infty(\widehat{G}) \mid \\
  \text{\((\psi\otimes\cnj{\chi}\bullet\cnj{\psi}) \cdot f\) is a
    \(\CONT_0\)-section of \(V\otimes \chi^*\cnj{V}\) on
    \(S\cap \cnj{\chi}S\)}\}.
\end{multline}
Here~\(\cnj{V}\) denotes the dual line bundle to~\(V\) and~\(\chi^*\)
pulls it back along the map \(\omega\mapsto \chi\omega\) to a line
bundle on \(\cnj{\chi}S\), so that the tensor product line bundle
\(V\otimes \chi^*\cnj{V}\) is defined on \(S\cap\cnj{\chi}S\); notice
also that \(\psi\otimes\cnj{\chi}\bullet\cnj{\psi}\) is a section of
this line bundle.  Since \(V\otimes \cnj{V}\) is a trivial line bundle
and \(\psi\otimes \cnj{\psi}\) is the constant section~\(1\),
\eqref{eq:Fell_fiber_at_chi} reduces to \eqref{eq:Fell_fiber_at_1} if
\(\chi=1\).

\begin{exmp}
  \label{exa:Ltwo_hatG_Zsquare_continued}
  We make \eqref{eq:Fell_fiber_at_chi} more explicit for the vector
  bundles~\(V_n\) over~\(\Torus^2\) considered in Example
  \ref{exa:Ltwo_hatG_Zsquare}.  It is convenient to describe
  \(\CONT(\Torus^2,V_n)\) as follows:
  \[
  \CONT(\Torus^2,V_n) =
  \{ f\in\CONT(\Torus\times\Rset) \mid
  f(z,t+1) = z^n f(z,t) \quad \forall\; z\in\Torus, t\in\Rset\}.
  \]
  The measurable section~\(\psi\) corresponds to \(\psi(z,t) \defeq z^{n
    \lfloor t\rfloor}\) in this picture, where \(\lfloor
  t\rfloor\in\Zset\) denotes the maximal integer \({}\le t\).  We fix
  \(\chi\in\Torus^2\) and represent it by \((a,b)\in\Torus\times\Rset\).
  Then
  \begin{multline*}
    \CONT(\Torus^2,\chi^*\cnj{V_n})
    = \{ f\in\CONT(\Torus\times\Rset) \mid
    f(a^{-1}z,t-b+1) = z^{-n} f(a^{-1}z,t-b)\}
    \\ = \{ f\in\CONT(\Torus\times\Rset) \mid
    f(z,t+1) = (az)^{-n} f(z,t) \quad \forall\; z\in\Torus, t\in\Rset\}
  \end{multline*}
  and hence
  \begin{multline*}
    \CONT(\Torus^2,V_n\otimes\chi^*\cnj{V_n})
    \\= \{f\in\CONT(\Torus\times\Rset)\mid
    f(z,t+1) = a^{-n} f(z,t) \quad\forall\; z\in\Torus, t\in\Rset\}.
  \end{multline*}
  The line bundle \(V_n\otimes\cnj{\chi}\bullet\cnj{V_n}\) must be
  trivial because~\(\Torus^2\) is connected, so that
  \(\chi^*\cnj{V_n}\cong \cnj{V_n}\).  The nowhere vanishing section
  \(\sigma_\chi(z,t)\defeq a^{-nt}\) provides an explicit trivialization
  for \(V_n\otimes\cnj{\chi}\bullet\cnj{V_n}\).  We can rewrite the
  condition \(f\cdot (\psi\otimes \cnj{\chi}\bullet\cnj{\psi}) \in
  \CONT(\Torus^2,V_n \otimes \chi^*\cnj{V_n})\) more explicitly as
  \[
  f\in (\psi\otimes \cnj{\chi}\bullet\cnj{\psi})^{-1}\cdot
  \CONT(\Torus^2,V_n\otimes\chi^*\cnj{V_n})
  = (\cnj{\psi}\otimes \cnj{\chi}\bullet\psi)\cdot\sigma_\chi \cdot
  \CONT(\Torus^2).
  \]
  We let
  \begin{multline*}
    h_\chi(z,t) \defeq
    (\cnj{\psi}\otimes \cnj{\chi}\bullet\psi)(z,t)\cdot \sigma_\chi(z,t)
    = z^{-n \lfloor t\rfloor}\cdot (az)^{n \lfloor t+b\rfloor}\cdot
    a^{-n t}
    \\= z^{n (\lfloor t+b\rfloor - \lfloor t\rfloor)}
    a^{n(\lfloor t+b\rfloor-t)}.
  \end{multline*}
  Notice that \(h_\chi(z,t+1)=h_\chi(z,t)\), so that~\(h_\chi\) defines
  a function on~\(\Torus^2\).  This function has jump discontinuities if
  \(t\in\Zset\) or \(t+b\in\Zset\) and is continuous otherwise.  Putting
  things together, we get
  \[
  \Bundle{B}_\chi(A,\Relc_{\Torus^2,V_n,\psi}) =
  \{M_f \circ M_{h_\chi} T_\chi^*\mid f\in \CONT(\Torus^2)\}.
  \]
\end{exmp}

\subsection{Spectrally Proper Algebras}
\label{sec:spec_prop}

The examples in \textsection\ref{sec:ex_Hilbert} show that a
\(G\)\nbd{}\Cstar{}algebra~\(A\) may have more than one or no continuous
spectral decomposition.  Now we consider a class of
\(G\)\nbd{}\Cstar{}algebras with a unique continuous spectral
decomposition.

Let~\(A\) be a \(G\)\nbd{}\Cstar{}algebra and let~\(\widehat{A}\) be its
\emph{spectrum}, that is, the set of unitary equivalence classes of
non-zero irreducible representations of~\(A\), endowed with the
hull-kernel topology.  The action~\(\alpha\) of~\(G\) on~\(A\) induces a
continuous action on~\(\widehat{A}\) by
\[
t\cdot\pi(a)\defeq \pi\circ \alpha_t^{-1}(a)
\qquad \forall\; t\in G, \pi\in \widehat{A}, a\in A.
\]
Let \(\Prim(A)\) be the \emph{primitive ideal space} of~\(A\) with the
Jacobson topology.  The group~\(G\) acts on \(\Prim(A)\) in an evident
fashion, and the canonical map \(\widehat{A}\to\Prim(A)\), \(\pi\mapsto
\ker\pi\), is \(G\)\nbd{}equivariant, continuous and open.  Even more,
it induces an isomorphism between the lattices of open subsets of
\(\widehat{A}\) and \(\Prim(A)\).

Let \(\Bundle{B} = (\Bundle{B}_\chi)_{\chi\in\widehat{G}}\) be a Fell
bundle over~\(\widehat{G}\).  Its
\emph{spectrum}~\(\widehat{\Bundle{B}}\) is its set of equivalence
classes of non-zero irreducible \Star{}representations.  The universal
property of \(A\defeq C^*(\Bundle{B})\) yields a natural homeomorphism
\(\widehat{\Bundle{B}} \cong \widehat{A}\) (see
\cite[VIII.17.3]{fell_doran}), which maps a representation \(\pi\colon
\Bundle{B}\to\Bound(\Hils)\) to its integrated form \(C^*(\Bundle{B})\to
\Bound(\Hils)\).  Naturality implies that the homeomorphism
\(\widehat{\Bundle{B}} \cong \widehat{A}\) is \(G\)\nbd{}equivariant for
the canonical \(G\)\nbd{}actions on both spaces.

The usual definition of proper group actions on locally compact
\emph{Hausdorff} spaces is extended in Definition 9.1 in~\cite{ralf1} to
group actions on locally \emph{quasi}-compact spaces such as
\(\Prim(A)\) and~\(\widehat{A}\).  Since it only uses the action on the
lattice of open subsets, \(\Prim(A)\) is a proper \(G\)\nbd{}space if
and only if~\(\widehat{A}\) is.

\begin{defn}[\cite{ralf1}]
  \label{def:sproper}
  A \(G\)\nbd{}\Cstar{}algebra~\(A\) is \emph{spectrally proper} if the
  \(G\)\nbd{}space \(\Prim(A)\) is proper.
\end{defn}

\begin{defn}
  \label{def:spec_prop_Fell}
  A Fell bundle~\(\Bundle{B}\) over~\(\widehat{G}\) is \emph{spectrally
    proper} if~\(\widehat{\Bundle{B}}\) is a proper \(G\)\nbd{}space (in the
  sense of Definition 9.1 in~\cite{ralf1}) or, equivalently, if the
  \(G\)\nbd{}\Cstar{}algebra \(C^*(\Bundle{B})\) is spectrally proper.
\end{defn}

\begin{thm}
  \label{the:sproper}
  The functor \(\Bundle{B}\mapsto C^*(\Bundle{B})\) is an equivalence of
  categories between the categories of spectrally proper Fell bundles
  over~\(\widehat{G}\) and of spectrally proper
  \(G\)\nbd{}\Cstar{}algebras.

  Hence spectrally proper \(G\)\nbd{}\Cstar{}algebra have a unique
  continuous spectral decomposition (up to canonical isomorphism).
  Spectrally proper Fell bundles are isomorphic as Fell bundles if and
  only if their cross-sectional \Cstar{}algebras are equivariantly
  isomorphic.
\end{thm}

\begin{proof}
  A spectrally proper \(G\)\nbd{}\Cstar{}algebra~\(A\) contains a unique
  relatively continuous, square-integrable, complete, dense
  subspace~\(\Relc_A\) by Theorem 9.1 in~\cite{ralf1}.  By our main
  theorem (Theorem \ref{the:main}), these subspaces correspond
  bijectively to continuous spectral decompositions of~\(A\).

  Since the subspace~\(\Relc_A\) is determined uniquely, any
  \(G\)\nbd{}equivariant \Star{}homomorphism \(A_1\to A_2\) maps
  \(\Relc_{A_1}\) to~\(\Relc_{A_2}\).  That is, the category of
  spectrally proper \(G\)\nbd{}\Cstar{}algebras is a full subcategory of
  the category of continuously square-integrable
  \(G\)\nbd{}\Cstar{}algebras.  Hence the equivalence of categories in
  Theorem \ref{the:main} yields the first assertion; it contains the
  remaining assertions.%\qed
\end{proof}

If~\(G\) is compact (and Abelian), then any \(G\)\nbd{}space is proper.
Hence Theorem \ref{the:sproper} specializes to the well-known assertion
that \(G\)\nbd{}\Cstar{}algebras for compact~\(G\) have a unique
spectral decomposition and that this provides an equivalence between the
categories of \(G\)\nbd{}\Cstar{}algebras and of Fell bundles
over~\(\widehat{G}\).

The unique continuous spectral decomposition of a spectrally proper
\(G\)\nbd{}\Cstar{}algebra is somewhat complicated to describe.
Therefore, we now restrict attention to \(G\)\nbd{}\Cstar{}algebras that
are \emph{proper in the sense of Kasparov} \cite{kasparov:novikov}, that
is, there is an essential \(G\)\nbd{}equivariant \Star{}homomorphism
from \(\CONT_0(X)\) to the center of \(\Mult(A)\) for some
\emph{Hausdorff} locally compact proper \(G\)\nbd{}space~\(X\).  It is
shown in \cite{Nilsen} that this is equivalent to the existence of a
continuous \(G\)\nbd{}map \(\Prim(A)\to X\) or to an isomorphism
\(A\cong\CONT_0(X,\Bundle{A})\) for some upper semi-continuous,
\(G\)\nbd{}equivariant field of \Cstar{}algebras \(\Bundle{A} =
(\Bundle{A_x})_{x\in X}\) over~\(X\).

Upper semi-continuity means that the functions \(X\ni x\mapsto
\norm{a_x}_{\Bundle{A_x}}\) are upper semi-continuous for all continuous
sections \((a_x)_{x\in X}\) of~\(\Bundle{A}\); in particular,
\(\norm{a_x}<\varepsilon\) for some \(x\in X\) implies that there is a
neighborhood~\(U\) of~\(x\) with \(\norm{a_y}<\varepsilon\) for all
\(y\in U\).  The \(G\)\nbd{}equivariance of the field of
\Cstar{}algebras means that the \(G\)\nbd{}action on~\(X\) lifts to
isomorphisms \(\alpha_t\colon A_x\to A_{t\cdot x}\) for all \(t\in G\),
\(x\in X\), that satisfy a suitable continuity condition.

Let \(\pi\colon \widehat{G}\times X\to X\) be the coordinate projection.
Pulling back the field~\(\Bundle{A}\) along~\(\pi\), we get a upper
semi-continuous field of \Cstar{}algebras \(\pi^*\Bundle{A}\) on
\(\widehat{G}\times X\); its \Cstar{}algebra of
\(\CONT_0\)\nbd{}sections is \(\CONT_0(\widehat{G},A)\), with the
obvious \(\CONT_0(\widehat{G}\times X)\)\brd{}module structure.

\begin{thm}
  \label{the:csd_proper_Kasparov}
  Let \(A=\CONT_0(X,\Bundle{A})\) for some upper semi-continuous
  \(G\)\nbd{}equivariant field of \Cstar{}algebras~\(\Bundle{A}\) over a
  proper \(G\)\nbd{}space~\(X\).  The fibers of the unique continuous
  spectral decomposition of~\(A\) are the spaces
  \(\Bundle{B}_\omega(A)\) of all bounded continuous sections \(b\in
  \CONT_b(X,\Bundle{A})\) that satisfy the following two conditions:
  \begin{itemize}
  \item \(\alpha_t\bigl(b(t^{-1}\cdot x)\bigr)= \braket{\omega}{t}\cdot
    b(x)\) for all \(t\in G\), \(x\in X\)

  \item \(G\backslash X\ni Gx\mapsto \norm{b(x)}\) belongs to
    \(\CONT_0(G\backslash X)\)
  \end{itemize}

  Let \(\pi^*\Bundle{A}\) be the pull-back of~\(\Bundle{A}\) along the
  projection \(\pi\colon \widehat{G}\times X\to X\).  A compactly
  supported section \(f\colon \widehat{G}\to \bigcup
  \Bundle{B}_\omega(A)\) is continuous if and only if the resulting
  section \((\omega,x)\mapsto f(\omega)(x)\) of \(\pi^*\Bundle{A}\) is
  (norm) continuous on \(\widehat{G}\times X\) and for each
  \(\varepsilon>0\) there is a compact subset \(C\subseteq X\) such that
  \(\norm{f(\omega)(x)} <\varepsilon\) for all \(\omega\in\widehat{G}\),
  \(x\in X\setminus G\cdot C\).

  The isomorphism \(C^*\bigl(\Bundle{B}(A)\bigr)\to A\) maps such a
  continuous section \((b_\omega)_{\omega\in\widehat{G}}\) to
  \(\int^\su_{\widehat{G}} b_\omega \dd{\omega}\).
\end{thm}

\begin{proof}
  Let \(A_c=\CONT_c(X,\Bundle{A})\) be the subspace of compactly
  supported continuous sections of~\(\Bundle{A}\).  Each element
  of~\(A_c\) is a linear combination of elements of the form \(a*a^*\)
  with \(a\in A_c\).  Since~\(A_c\) is square-integrable and
  relatively continuous, we conclude that~\(A_c\) is a relatively
  continuous subset of integrable elements.  Furthermore, \(A_c\) is a
  \(G\)\nbd{}invariant ideal in~\(A\).  Therefore, the fibers of the
  unique continuous spectral decomposition of~\(A\) are
  \(\Bundle{D}_\omega\defeq \clos{E_\omega(A_c)}\), and the continuous
  sections are generated by \(\omega\mapsto E_\omega(a)\) for \(a\in
  A_c\) as in Lemma \ref{lem:Fell_cont_section}.  We must compare this
  with the fibers~\(\Bundle{B}_\omega\) and spaces of continuous
  sections described in the statement of the theorem.

  If \(a\in A_c\) is supported in~\(C\), then \(b_\omega\defeq
  E_\omega(a)\) is supported in \(G\cdot C\) for all
  \(\omega\in\hat{G}\).  Hence \(\norm{b_\omega(x)}\) vanishes for
  \(G\cdot x\to\infty\) in \(G\backslash X\), even uniformly
  in~\(\omega\).  Furthermore, \(b_\omega\) is a continuous section
  of~\(\Bundle{A}\) because the integrand \(\alpha_t\bigl(a(t^{-1}\cdot
  x)\bigr)\cdot \cnj{\braket{\omega}{t}}\) is compactly supported
  uniformly for \(x\in U\) for any relatively compact subset
  \(U\subseteq X\).

  It is also clear that \(\alpha_t\bigl(b_\omega(t^{-1}\cdot x)\bigr) =
  \braket{\omega}{t}\cdot b_\omega(x)\) for all \(t\in G\), \(x\in X\),
  \(\omega\in\widehat{G}\): this is equivalent to
  \(b_\omega\in\Mult_\omega(A)\) and follows from
  \eqref{eq:x-spectral_sub}.  This yields \(E_\omega(A_c)\subseteq
  \Bundle{B}_\omega\).  We even get
  \(\Bundle{D}_\omega\subseteq\Bundle{B}_\omega\)
  because~\(\Bundle{B}_\omega\) is a closed linear subspace of
  \(\Mult_\omega(A)\).  Similar arguments show that continuous sections
  of~\(\Bundle{D}_\omega\) satisfy the conditions in the statement of
  the theorem.

  Conversely, we must find \(a\in A_c\) with \(\norm{b-
    E_\omega(a)}<\varepsilon\) for \(b\in\Bundle{B}_\omega\) and
  \(\varepsilon>0\).  For \(h\in\CONT_c(X)\), we view
  \(E_1(h)\in\Mult\bigl(\CONT_0(X)\bigr) \cong \CONT_b(X)\) as a
  function on~\(X\):
  \[
  E_1(h)(x) = \int_G h(t^{-1}\cdot x)\dd{t}.
  \]
  The functions \(E_1(h)\) with \(h\in\CONT_c(X)\) form an ideal in
  \(\CONT_0(G\backslash X)\), which is dense by the Stone--Weierstraß
  Theorem.  Thus we get an approximate unit for \(\CONT_0(G\backslash
  X)\) of the form \(\bigl(E_1(h_i)\bigr)\).  Since
  \(b\in\Mult_\omega(A)\), we have
  \[
  E_\omega(h_i\cdot b) = E_1(h_i)\cdot b.
  \]
  We get \(\norm{E_\omega(h_i\cdot b)-b}<\varepsilon\) for sufficiently
  large~\(i\) because \(E_1(h_i)\) is an approximate identity and the
  norm of~\(b\) is in \(\CONT_0(G\backslash X)\).  Since \(h_i\cdot b\in
  A_c\), this finishes the proof that
  \(\Bundle{D}_\omega=\Bundle{B}_\omega\) for all
  \(\omega\in\widehat{G}\).

  Now we let \((b_\omega)_{\omega\in\widehat{G}}\) be a compactly
  supported continuous section of~\(\Bundle{B}_\omega\) as described in
  the theorem.  The uniform vanishing at infinity in the direction of
  \(G\backslash X\) allows us to approximate~\(b_\omega\) uniformly by a
  compactly supported function on \(\widehat{G}\times X\).  Since the
  space of continuous sections of~\(\Bundle{D}\) is a Banach
  \(\CONT_0(\widehat{G}\times X)\)\brd{}module, \((b_\omega)\) is a
  continuous section of~\(\Bundle{D}\) if and only if we can approximate
  it \emph{locally} by functions of the form \(E_\omega(a)\) for \(a\in
  A_c\); even more, the upper semi-continuity of the norm shows that it
  suffices to approximate it \emph{pointwise}.  This reduces us to the
  case of fixed \(\omega\in\widehat{G}\), which we have treated above.
  Hence the continuous sections of the Fell bundle~\(\Bundle{D}\) are as
  described in the statement of the theorem.

  Finally, the formula for the isomorphism \(C^*(\Bundle{B})\to A\) is a
  special case of Lemma~\ref{lem:iso_Fell}.
\end{proof}

\subsection{Proper actions and commutative Fell bundles}
\label{sec:commutative_Fell}

Now we consider an even more special case.  A Fell bundle
over~\(\widehat{G}\) is called \emph{commutative} if \(b_\chi\cdot b_\omega
= b_\omega\cdot b_\chi\) for all \(b_\chi\in\Bundle{B}_\chi\),
\(b_\omega\in\Bundle{B}_\omega\).  Equivalently, \(C^*(\Bundle{B})\) is
commutative.  In this case, the spectrum~\(\widehat{\Bundle{B}}\) is a
Hausdorff locally compact space and agrees with the space of characters
of~\(\Bundle{B}\) or, equivalently, of \(C^*(\Bundle{B})\).

\begin{thm}
  \label{the:continuous_Fell}
  Isomorphism classes of commutative Fell bundles over~\(\widehat{G}\)
  correspond bijectively to isomorphism classes of proper
  \(G\)\nbd{}spaces.
\end{thm}

\begin{proof}
  The Gelfand--Naimark Theorem yields \(C^*(\Bundle{B}) \cong
  \CONT_0(\widehat{\Bundle{B}})\).  Proposition \ref{pro:CONT_c_relc} shows
  that the \(G\)\nbd{}action on \(\CONT_0(\widehat{\Bundle{B}})\) is
  square-integrable.  It was already observed by Marc Rieffel
  \cite{rieffel2} that this implies that~\(\widehat{\Bundle{B}}\) is a
  proper \(G\)\nbd{}space.  As a result, all commutative Fell bundles
  are spectrally proper.  Now Theorem \ref{the:sproper} shows that
  \(\Bundle{B}\mapsto C^*(\Bundle{B})\) is an equivalence of categories
  between the category of commutative Fell bundles and the category of
  commutative, spectrally proper \(G\)\nbd{}\Cstar{}algebras; here we
  use that a continuous spectral decomposition of a commutative
  \(G\)\nbd{}\Cstar{}algebra is again commutative.  Of course, a
  commutative \emph{spectrally proper} \(G\)\nbd{}\Cstar{}algebra is of
  the form \(\CONT_0(X)\) for a \emph{proper} locally compact
  \(G\)\nbd{}space~\(X\), and \(X\) and \(\CONT_0(X)\) determine each
  other uniquely up to canonical isomorphism.%\qed
\end{proof}

We warn the reader that the map \(\CONT_0(X)\mapsto X\) from commutative
\Cstar{}algebras to locally compact spaces is not functorial; for
instance, the zero homomorphism \(\CONT_0(Y)\to\CONT_0(X)\) does not
correspond to a map \(X\to Y\).  Instead, \Star{}homomorphisms
\(\CONT_0(Y)\to\CONT_0(X)\) correspond bijectively to continuous pointed
maps \(X_+\to Y_+\), where \(X_+\) and~\(Y_+\) denote the one-point
compactifications of \(X\) and~\(Y\).

Theorem \ref{the:continuous_Fell} improves upon Corollary 12.5
in~\cite{exel1}.  It is also related to Theorem X.5.13
in~\cite{fell_doran}, which asserts that \(\Bundle{B}\mapsto
\widehat{\Bundle{B}}\) is a bijection between isomorphism classes of
\emph{saturated} commutative Fell bundles over~\(\widehat{G}\) and locally
compact principal \(G\)\nbd{}bundles or, equivalently, free and proper
\(G\)\nbd{}spaces; saturated means that \(\cspn
{}(\Bundle{B}_\chi\cdot\Bundle{B}_\eta)=\Bundle{B}_{\chi\eta}\), for all
\(\chi,\eta\in\widehat{G}\).

\begin{prop}
  \label{pro:Fell_bundle_proper}
  Let~\(X\) be a proper \(G\)\nbd{}space.  The unique continuous
  spectral decomposition of \(\CONT_0(X)\) has as its fibers the
  spaces~\(\Bundle{B}_\omega\) of all \(f\in\CONT_b(X)\) that satisfy
  the following two conditions:
  \begin{itemize}
  \item \(\alpha_t(f)=\braket{\omega}{t} \cdot f\) for all \(t\in G\),
    \(\omega\in\widehat{G}\)

  \item \(\abs{f}\in\CONT_0(G\backslash X)\)

  \end{itemize}
  A compactly supported section of this Fell bundle is continuous if and
  only if it is continuous as a function \(f\colon \widehat{G}\times
  X\to\Cset\) and if, for each \(\varepsilon>0\), there is \(C\subseteq
  X\) with \(\abs{f(x)}<\varepsilon\) for \(x\in X\setminus G\cdot C\).
\end{prop}

\begin{proof}
  This is the special case of Theorem~\ref{the:csd_proper_Kasparov}
  where~\(\Bundle{A}\) is the trivial bundle with fiber~\(\Cset\).
\end{proof}


\begin{thebibliography}{9}
\bibitem{exel2} R. Exel,
  Unconditional Integrability for Dual Actions,
  \textit{Bol. Soc. Brasil. Mat. (N.S.)} \textbf{30} (1999), 99--124.

\bibitem{exel1} R. Exel,
  Morita--Rieffel Equivalence and Spectral Theory for Integrable
  Automorphism Groups of C*-algebras,
  \textit{J. Funct. Anal.} \textbf{172} (2000), 404--465.

\bibitem{fell_doran} J. M. G. Fell, R. S. Doran,
  ``Representations of \(*\)-Algebras, Locally Compact Groups, and
    Banach \(*\)-Algebraic Bundles,''
  Academic Press, 1988.

\bibitem{kasparov:novikov} G. Kasparov,
  Equivariant $KK$-theory and the Novikov conjecture,
  \textit{Invent. Math.} \textbf{91} (1988), 147--201.

\bibitem{ralf2} R. Meyer,
  Equivariant Kasparov Theory and Generalized Homomorphisms,
  \textit{K\nobreakdash-Theory} \textbf{21} (2000), 201--228.

\bibitem{ralf1} R. Meyer,
  Generalized Fixed Point Algebras and Square-Integrable Groups Actions,
  \textit{J. Funct. Anal.} \textbf{186} (2001), 167--195.

\bibitem{Nilsen} M. Nilsen,
  \(C^*\)-bundles and \(C_0(X)\)-algebras,
  \textit{Indiana Univ. Math. J.} \textbf{45} (1996), 463--477.

\bibitem{rieffel2} M. A. Rieffel,
  Proper actions of Groups on \Cstar{}algebras, \textit{in}
  ``Mappings of Operator Algebras (Philadelphia, PA, 1988),''
  141--182, Birkh\"auser, 1990.

\bibitem{rieffel1} M. A. Rieffel,
  Integrable and Proper Actions on \Cstar{}algebras, and
    Square-Integrable Representations of Groups,
  math.OA/9809098 v2, 22 June, 1999.

\end{thebibliography}
\end{document}